\documentclass[10pt]{article}
\usepackage{latexsym}
\usepackage{amsfonts}
\usepackage{amssymb}
\usepackage[utf8]{inputenc}
\usepackage[active]{srcltx} %SRC Specials for DVI Searching
\usepackage{color}
\usepackage[all]{xy}
\usepackage{amsthm}

\usepackage{amsmath}
\newcommand{\Ima}{\text{Im}}

\newcommand{\diver}{{\rm div}}

\renewcommand{\span}{{\rm span}}

\newcommand{\ba}{\begin{array}}\newcommand{\ea}{\end{array}}

\newcommand{\ns}{\rm}

\newcommand{\Id}{\mathrm{Id}}

\newcommand{\nse}{\kern-3pt\ns$=$}\newcommand{\qd}{\hfill$\Box$\medbreak}

\newcommand{\ext}{\raise1pt\hbox{$\ts\bigwedge$}}

\newcommand{\ts}{\textstyle}
\newcommand{\rf}[1]{(\ref{#1})}
\newcommand{\chii}{\raise2pt\hbox{$\chi$}}

\newcommand{\Fg}{\mbox{${\cal F}\kern-2pt_g$}}

\newcommand{\Mg}{\mbox{${\cal M}\kern-2pt_g$}}
\newcommand{\Ng}{\mbox{${\cal N}\kern-2pt_g$}}
\newcommand{\V}{V\kern-1pt}

\newcommand{\Gg}{\mbox{${\cal G}\kern-2pt_g$}}
\newcommand{\tr}{\mbox{\bf tr}}

\newcommand{\Res}[2]{\hbox{\ns Res}\kern-16pt\lower5pt\hbox{\footnotesize$_{#1}$}\kern2pt\left[#2\right]}
\newcommand{\qk}{quaternion-K\"ahler\kern2pt}\renewcommand{\,}{\kern1pt}

\newcommand{\End}{{\rm End}}

\newcommand{\dirac}{/\kern-5pt\partial}

\newcommand{\Spinc}{Spin\hbox{$^c$}\,\,}

\newcommand{\lra}{\longrightarrow}
\renewcommand{\ts}{\textstyle}
0
\newtheorem{theo}{Theorem}[section]
\newtheorem{defi}{Definition}[section]\newtheorem{lemma}{Lemma}[section]

\newtheorem{corol}{Corollary}[section]
\newtheorem{prop}{Proposition}[section]\def\frac#1#2{{#1\over#2}}
\def\be#1\ee{\begin{equation}#1\end{equation}}

\allowdisplaybreaks

\begin{document}

\title{Spinorially twisted Spin structures, III: CR structures}

\author{
Rafael Herrera\footnote{Centro de
Investigaci\'on en Matem\'aticas, A. P. 402,
Guanajuato, Gto., C.P. 36000, M\'exico. E-mail: rherrera@cimat.mx}
\footnote{Partially supported by 
grants of CONACyT, LAISLA (CONACyT-CNRS), and the IMU Berlin Einstein Foundation
Program},
Roger Nakad\footnote{Notre Dame University-Louaiz\'e, Faculty of Natural and Applied Sciences, Department of Mathematics and Statistics, P.O. Box 72, Zouk Mikael, Lebanon. 
E-mail: rnakad@ndu.edu.lb}
\footnote{Partially supported by a Fellowship of the IMU Berlin Einstein Foundation
Program} \,\, and
Iv\'an T\'ellez\footnote{Centro de
Investigaci\'on en Matem\'aticas, A. P. 402,
Guanajuato, Gto., C.P. 36000, M\'exico. E-mail: tellezito@cimat.mx}
\footnote{Supported by 
a CONACyT scholarship}
 }

\date{}

\maketitle

\vspace{-20pt}

{
{{\abstract{
We develop a spinorial description of CR structures of arbitrary codimension.
More precisely, we characterize almost CR structures of arbitrary codimension on (Riemannian) manifolds by the existence of a Spin$^{c, r}$ structure 
carrying a partially pure spinor field. We study various integrability conditions of the almost CR structure in our spinorial setup, including the classical integrability of a CR structure 
as well as those implied by Killing-type conditions on the partially pure spinor field.  In the codimension one case, 
we develop a spinorial description of strictly pseudoconvex CR manifolds, metric contact manifolds and  Sasakian manifolds. 
Finally, we study hypersurfaces of K\"ahler manifolds via partially pure Spin$^c$ spinors.
}
}

}}

\section{Introduction}

Spinors  have played an important role in both physics and mathematics 
ever since they were discovered by \'E. Cartan in $1913$. 
We refer the reader to Hitchin's seminal paper \cite{hitchin}, Friedrich's textbook \cite{Friedrich}, as well as to 
\cite{witten, morgan} for the more recent development of Seiberg-Witten theory and its notorious
results on $4$-manifold geometry and topology.

The starting point of this paper was our interest in characterizing and studying 
CR structures (of arbitrary codimension) by means of twisted Spin structures and spinors.
Our main motivation has been the relation of almost complex structures with
``classical pure spinors'' (and Spin$^c$ structures). 
Cartan defined {\em pure spinors} \cite{cartan,
cartan1, cheva} in order to characterize (almost) complex structures and, almost one hundred years
later, they are still being used in related geometrical problems \cite{borisov}. 
Furthermore, these spinor fields
have been related to the notion of calibrations on a Spin manifold by Harvey and Lawson
\cite{har, dadok}, since distinguished differential forms are naturally 
associated to a spinor field and, in particular, give rise
to special differential forms on immersed hypersurfaces. Pure spinors are also present in the
Penrose formalism 
in General Relativity as they are implicit in Penrose's notion of ``flag planes'' \cite{pen1,
pen2, pen3}.

The notion of abstract CR structures in odd dimensions generalizes that of 
complex structure in even dimensions. This notion aims to describe intrinsically the
property of being a hypersurface of a complex space form. This is done by distinguishing a
distribution whose sections play the role of the holomorphic vector fields tangent to the
hypersurface. 
There exists also the notion of almost CR structure of arbitrary codimension, in which a fixed codimension subbundle
of the tangent bundle carries a complex structure.
It has been proved that
every codimension one, strictly
pseudoconvex CR manifold has a canonical Spin$^c$ structure \cite{petit}. 
Naturally, this led us to ask if it is possible to characterize almost CR structures of arbitrary codimension (and a choice of compatible metric)
by means of a twisted Spin structure carrying a special spinor field.  

We developed the algebraic background of twisted partially pure spinors in \cite{Herrera-Tellez} 
which we recall briefly in the second section.
Let us recall, in particular, the definition of twisted Spin group 
\[Spin^{c,r}(n)={Spin(n)\times Spin^c(r)\over \{\pm (1,1)\}},\]
which will be the structure group for the twisted Spin structures (cf. Definition \ref{defi-twisted-spin-structure}), and whose representations
contain the partially pure spinors. Note that $r$ will eventually be the codimension of an almost CR structure.
Such twisted Spin structures involve not only the principal bundle of orthonormal frames, 
but also two auxiliary principal bundles. The need for such structures stems from the fact that
there are manifolds which are neither Spin nor Spin$^c$. Subsection \ref{subsec: existence of non-reducible twisted Spin structures} 
is devoted to showing that there are triples of principal bundles admitting Spin$^{c,r}$ structures.

The existence of a partially pure spinor field $\phi$ on a Riemannian Spin$^{c,r}$ manifold $M^n$ implies the
splitting of the tangent bundle $TM$ into two orthogonal distributions $V^\phi$ and $(V^\phi)^\perp$, where the former
is endowed with an automorphism $J^\phi$ satisfying $(J^\phi)^2 = -{\rm Id}_{V^\phi}$, i.e. $M$ has an almost CR hermitian
structure. In fact, the converse is also
true (cf. Theorem \ref{theo: characterization almost CR}). 
Furthermore, we characterize the integrability condition of a CR structure (with metric) by
an equation involving covariant derivatives of the partially pure spinor (cf. Theorem \ref{theo: integrable}).
We proceed to study other natural ``integrability conditions'' of the partially pure spinor field,
such as being parallel in the $V^\phi$ directions (cf. Theorem \ref{theo: D-parallel}), or 
being Killing in the $(V^\phi)^\perp$ directions (cf. Theorem \ref{theo: Dperp Killing}), etc.
We present a family of homogeneous spaces as examples for the different theorems.

As mentioned before, the relevant group for codimension one almost CR structures is Spin$^{c,1}(n)=$ Spin$^c(n)$.
Thus, we prove that partially pure
spinors  appear
 naturally and implicitly in extrinsic Spin$^c$ geometry: consider a K\"ahler manifold endowed with
a Spin$^c$  structure carrying a parallel spinor $\psi$.  It is known that the restriction $\phi $ of the
parallel spinor $\psi$ to a real oriented hypersurface $M$   satisfies 
\[
 \nabla^M_X\phi = -{1\over 2}\, {\rm II}(X)\bullet\phi, 
\]
where ${\rm II}$ denotes the second fundamental form of $M$, $\nabla^M$ is the Spin$^c$ covariant
derivative on $M$ and ``$\bullet$'' the Clifford multiplication on $M$ \cite{morel, r2}. Moreover,  the spinor $\phi$  is partially pure and integrable (see
Theorem \ref{kaa}).

The paper is organized as follows. In Section \ref{sec:preliminaries}, we recall the background material for 
the definition of partially pure spinors, and describe the isotropy representation of a family of homogeneous spaces 
(partial flag manifolds) that will be used throughout the paper.
In Section \ref{sec: doubly twisted Spin structures}, 
we define Spin$^{c,r}$ structures, study their existence, define twisted Dirac and Laplacian operators, prove 
some curvature identities and a 
Schr\"odinger-Lichnerowicz-type formula, and derive some Bochner-type results. 
In Section \ref{sec: CR structures}, we give the spinorial characterization of (almost) CR hermitian structures, and 
examine various integrability conditions and their geometrical consequences.
In Section \ref{sec: codimension one CR structures}, we return to the codimension one case and examine in our spinorial context
(strictly) pseudoconvex CR manifolds, metric contact manifolds and Sasakian manifolds, and explore
extrinsic geometry questions including immersion theorems (Theorems \ref{kaa} and
\ref{immersion}).

{\bf Acknowledgments}. 
The authors are grateful to Oussama Hijazi for his encouragement and valuable
comments. The authors thank  Helga Baum and the Institute of Mathematics of the University of
Humboldt-Berlin  for their hospitality and support. 
The first author would also like to thank the hospitality and support of the 
International Centre for Theoretical Physics and the  Institut des Hautes \'Etudes Scientifiques. 
The second author gratefully acknowledges  the support and hospitality of the 
Centro de Investigaci\'{o}n en Matem\'{a}ticas A.C. (CIMAT).

\section{Preliminaries}\label{sec:preliminaries}

In this section, we briefly recall basic facts about  Clifford algebras, the Spin
group and the standard Spin representation \cite{Friedrich}. We also define the twisted Spin groups and
representations, the antisymmetric 2-forms and endomorphisms associated to a twisted spinor, 
recall the definition of twisted partially pure spinor, and describe the isotropy representations of certain 
homogeneous spaces that will furnish examples later on.

\subsection{Clifford algebras, the Spin group and representation}

Let $Cl_n$ denote the Clifford algebra generated by the orthonormal vectors
$e_1, e_2, \ldots, e_n\in \mathbb{R}^n$ 
subject to the relations
\begin{eqnarray*}
e_j e_k + e_k e_j &=& -2\left< e_j,e_k\right>,
\end{eqnarray*}
where $\big< , \big>$ denotes the standard inner product in $\mathbb{R}^n$.
Let
\[\mathbb{C}l_n=Cl_n\otimes_{\mathbb{R}}\mathbb{C}\]
denote the complexification of $Cl_n$. The Clifford algebras are isomorphic to matrix algebras
\[\mathbb{C}l_n\cong \left\{
                     \begin{array}{ll}
                     \End(\mathbb{C}^{2^k}), & \mbox{if $n=2k$,}\\
                     \End(\mathbb{C}^{2^k})\oplus\End(\mathbb{C}^{2^k}), & \mbox{if $n=2k+1.$}
                     \end{array}
\right.
\]
The map
\[\kappa:\mathbb{C}l_n \lra \End(\Delta_n)\]
is defined to be either the above mentioned isomorphism if $n$ is even, or the isomorphism followed by
the projection onto the first summand if $n$ is odd. 
An expression for $\kappa$ can be given explicitly using the matrices
\[\Id = \left(\begin{array}{ll}
1 & 0\\
0 & 1
      \end{array}\right),\quad
g_1 = \left(\begin{array}{ll}
i & 0\\
0 & -i
      \end{array}\right),\quad
g_2 = \left(\begin{array}{ll}
0 & i\\
i & 0
      \end{array}\right),\quad
T = \left(\begin{array}{ll}
0 & -i\\
i & 0
      \end{array}\right).
\]
In terms of the generators $e_1, \ldots, e_n$, $\kappa$ is given by
\begin{eqnarray}
e_1&\mapsto& \Id\otimes \Id\otimes \ldots\otimes \Id\otimes \Id\otimes g_1,\nonumber\\
e_2&\mapsto& \Id\otimes \Id\otimes \ldots\otimes \Id\otimes \Id\otimes g_2,\nonumber\\
e_3&\mapsto& \Id\otimes \Id\otimes \ldots\otimes \Id\otimes g_1\otimes T,\nonumber\\
e_4&\mapsto& \Id\otimes \Id\otimes \ldots\otimes \Id\otimes g_2\otimes T,\nonumber\\
\vdots && \dots\nonumber\\
e_{2k-1}&\mapsto& g_1\otimes T\otimes \ldots\otimes T\otimes T\otimes T,\nonumber\\
e_{2k}&\mapsto& g_2\otimes T\otimes\ldots\otimes T\otimes T\otimes T,\nonumber
\end{eqnarray}
and, if $n=2k+1$, 
\[ e_{2k+1}\mapsto i\,\, T\otimes T\otimes\ldots\otimes T\otimes T\otimes T.\]
The vectors 
\[u_{+1}={1\over \sqrt{2}}(1,-i)\quad\quad\mbox{and}\quad\quad u_{-1}={1\over \sqrt{2}}(1,i),\]
form a unitary basis of $\mathbb{C}^2$ with respect to the standard Hermitian product.
Thus
\[\{u_{\varepsilon_1,\ldots,\varepsilon_k}=u_{\varepsilon_1}\otimes\ldots\otimes
u_{\varepsilon_k}\,\,|\,\, \varepsilon_j=\pm 1,
j=1,\ldots,k\},\]
is a unitary basis of $\Delta_n=\mathbb{C}^{2^k}$
with respect to the naturally induced Hermitian product.

The Clifford multiplication is defined by
\begin{eqnarray*}
\mu_n:\mathbb{R}^n\otimes \Delta_n &\lra&\Delta_n\\ 
x \otimes \psi &\mapsto& \mu_n(x\otimes \psi)=x\cdot\psi :=\kappa(x)(\psi).
\end{eqnarray*}
Aditionally, the maps 
\[\alpha\left(\begin{array}{c}
z_1\\
z_2
              \end{array}
\right) = \left(\begin{array}{c}
-\overline{z}_2\\
\overline{z}_1
              \end{array}\right), \qquad \beta\left(\begin{array}{c}
z_1\\
z_2
              \end{array}
\right) = \left(\begin{array}{c}
\overline{z}_1\\
\overline{z}_2
              \end{array}\right),\]
define quaternionic and real structures, respectively, on $\mathbb{C}^2$. Using $\alpha$ and $\beta$, real or quaternionic structures $\gamma_n$ are built on $\Delta_n=(\mathbb{C}^2)^{\otimes
[n/2]}$, for $n\geq 2$, as follows
\[
\begin{array}{cclll}
 \gamma_n &=& (\alpha\otimes\beta)^{\otimes 2k} &\mbox{if $n=8k,8k+1$}& \mbox{(real),} \\
 \gamma_n &=& \alpha\otimes(\beta\otimes\alpha)^{\otimes 2k} &\mbox{if $n=8k+2,8k+3$}&
\mbox{(quaternionic),} \\
 \gamma_n &=& (\alpha\otimes\beta)^{\otimes 2k+1} &\mbox{if $n=8k+4,8k+5$}&\mbox{(quaternionic),} \\
 \gamma_n &=& \alpha\otimes(\beta\otimes\alpha)^{\otimes 2k+1} &\mbox{if $n=8k+6,8k+7$}&\mbox{(real).}
\end{array}
\]
The Spin group $Spin(n)\subset Cl_n$ is the subset 
\[Spin(n) =\{x_1x_2\cdots x_{2l-1}x_{2l}\,\,|\,\,x_j\in\mathbb{R}^n, \,\,
|x_j|=1,\,\,l\in\mathbb{N}\},\]
endowed with the product of the Clifford algebra.
It is a Lie group and its Lie algebra is
\[\mathfrak{spin}(n)=\mbox{span}\{e_ie_j\,\,|\,\,1\leq i< j \leq n\}.\]
Recall that the Spin group $Spin(n)$ is the universal double cover of $SO(n)$, $n\ge 3$. For $n=2$
we consider $Spin(2)$ to be the connected double cover of $SO(2)$.
The covering map will be denoted by 
\[\lambda_n:Spin(n)\rightarrow SO(n).\]
Its differential is given
by $(\lambda_n)_{*}(e_ie_j) = 2E_{ij}$, where $E_{ij}=e_i^*\otimes e_j - e_j^*\otimes e_i$ is the
standard basis of the skew-symmetric matrices and $e^*$ denotes the metric dual of the vector $e$.
Furthermore, we will abuse the notation and also denote by $\lambda_n$ the induced representation on
$\ext^*\mathbb{R}^n$.

The restriction of $\kappa$ to $Spin(n)$ defines the Lie group representation
\[
\kappa_n:Spin(n)\lra GL(\Delta_n),\]
which is special unitary. We have the corresponding Lie algebra representation
\[
\kappa_{n*}:\mathfrak{spin}(n)\lra \mathfrak{gl}(\Delta_n).\] 

{\bf Remark}.
For the sake of notation we will set
\[SO(0)=\{1\},\quad\quad SO(1)=\{1\},\]
\[Spin(0)=\{\pm1\},\quad\quad Spin(1)=\{\pm1\},\]
and
\[\Delta_0 = \Delta_1 =\mathbb{C}\]
a trivial $1$-dimensional representation.

The Clifford multiplication $\mu_n$ 
% has the following properties:
% \begin{itemize}
%  \item It 
is skew-symmetric with respect to the Hermitian product
\[\left<x\cdot\psi_1 , \psi_2\right> 
=-\left<\psi_1 , x\cdot \psi_2\right>,
\] 
is $Spin(n)$-equivariant
% map of $Spin(n)$ representations.
%  \item $\mu_n$ 
and  can be extended to a $Spin(n)$-equivariant map 
\begin{eqnarray*}
\mu_n:\ext^*(\mathbb{R}^n)\otimes \Delta_n &\lra&\Delta_n\\ 
\omega \otimes \psi &\mapsto& \omega\cdot\psi.
\end{eqnarray*}
% of $Spin(n)$ representations. 
% \end{itemize}

\subsection{Twisted Spin groups}

Consider the following groups:
\begin{itemize}
 \item By using the unit complex numbers $U(1)$,
the Spin group can be twisted \cite{Friedrich}
\[Spin^c(n) =  (Spin(n) \times U(1))/\{\pm (1,1)\} =
Spin(n) \times_{\mathbb{Z}_2} U(1),\]
with Lie algebra
\[\mathfrak{spin}^c(n)=\mathfrak{spin}(n)\oplus i\mathbb{R}.\]

\item By using $Spin^c(r)$ define
\begin{eqnarray*}
Spin^{c,r}(n) &=&  (Spin(n) \times Spin^c(r))/\{\pm (1,1)\} \\
&=& Spin(n) \times_{\mathbb{Z}_2} Spin^c(r),
\end{eqnarray*}
where $r\in\mathbb{N}$, whose
Lie algebra is
\[\mathfrak{spin}^c(n)=\mathfrak{spin}(n)\oplus \mathfrak{spin}(r)\oplus i\mathbb{R}.\]
It fits into the exact sequence
\[1\lra \mathbb{Z}_2\lra Spin^{c,r}(n)\xrightarrow{\lambda_{n,r,2}} SO(n)\times SO(r)
\times U(1)\lra 1,\]
where
\begin{eqnarray*}
(\lambda_{n,r,2})[g,[h,z]]  &=& (\lambda_n(g),\lambda_r(h),z^2).
\end{eqnarray*}
\item Let 
$\widehat{Spin^c(r)}$ denote the standard copy of $Spin^c(r)$  in $Spin^{c,r}(n)$ given by elements
of the form
$[1,[h,z]]$ where $h\in Spin(r)$ and $z\in U(1)$.
\end{itemize}

\vspace{.1in} 

{\bf Remark}. For $r=0,1$,  $Spin^{c,r}(n)=Spin^c(n)$.

\begin{lemma}\label{lemma:subgroup2} {\rm \cite{Herrera-Tellez}}
Let $r\in \mathbb{N}$.
There exists a monomorphism $h:U(m)\times SO(r) \hookrightarrow
Spin^{c,r}(2m+r)$
such that the following diagram commutes
\[
\xymatrix{
  & Spin^{c,r}(2m+r) \ar[d]\\
U(m)\times SO(r) \ar[ur] \ar[r] & SO(2m+r)\times SO(r)\times U(1)
}
\]
% \[
% \xymatrix{
% \pi_1(G)\ar[d]_h \ar[r]^-{i_{\#}} & \pi_1(SO(n))\\
% \pi_1(Q) \ar[ur]_{f}
% }
% \]
% 
\end{lemma}
\qd

% \begin{lemma}\label{factorization} {\rm \cite{Herrera-Tellez}}
% Let $r\in\mathbb{N}$.
%  The standard representation $\Delta_{2m+r}$ of Spin$(2m+r)$ decomposes as follows
% \[\Delta_{2m+r} = \Delta_r\otimes\Delta_{2m}^+ \,\,\oplus\,\,\Delta_r\otimes\Delta_{2m}^-,\]
% with respect to the subgroup
% $Spin(2m)\times_{\mathbb{Z}_2} Spin(r) \subset Spin(2m+r)$. 
% \end{lemma}
% \qd

\begin{lemma}\label{lemma: fundamental group SpinCR} Let $n\geq 3$.
\begin{itemize}
 \item For $r\geq 3$
\[\pi_1(Spin^{c,r}(n))=\mathbb{Z}_2\oplus\mathbb{Z}.\] 
 \item For $r=2$
\[\pi_1(Spin^{c,2}(n))=\mathbb{Z}\oplus\mathbb{Z}.\] 
\end{itemize}
\end{lemma}
{\em Proof}. 
For $r\geq 3$, consider the universal cover
\[\begin{array}{c}
Spin(n)\times Spin(r)\times \mathbb{R}\\
\downarrow\\
Spin^{c,r}(n)
  \end{array}
\]
The preimage of $[1,[1,1]]\in$ Spin$^{c,r}(n)$ is
\[ \left<(-1,-1,0),(1,-1,1)\right>\subset Spin(n)\times Spin(r)\times \mathbb{R}.\]

For $r=2$, consider the universal cover
\[\begin{array}{c}
Spin(n)\times \mathbb{R}\times \mathbb{R}\\
\downarrow\\
Spin^{c,2}(n)
  \end{array}
\]
The preimage of $[1,[1,1]]\in Spin^{c,2}(n)$ is
\[ \left<(-1,1,0),(-1,0,1)\right>\subset Spin(n)\times \mathbb{R}\times \mathbb{R}.\]
\qd

\subsection{Twisted Spin representations}

Consider the following twisted representations:
\begin{itemize}
 \item The Spin representation $\Delta_n$ extends to a representation of Spin$^c(n)$ by letting
\begin{eqnarray*}
Spin^c(n)&\longrightarrow& GL(\Delta_n)\\
\,[g,z]  &\mapsto& z\kappa_n(g)=:zg.
\end{eqnarray*}

\item The twisted Spin$^{c,r}(n)$ representation is given by
\begin{eqnarray*}
\kappa_n^{c,r}:Spin^{c,r}(n)&\longrightarrow& GL(\Delta_r\otimes
\Delta_n)\\
\,[g,[h,z]]  &\mapsto& z \, \kappa_r(h)\otimes \kappa_n(g)=:zh\otimes g
\end{eqnarray*}
which is also unitary with respect to the natural Hermitian metric.
\item For $r=0,1$, the twisted Spin representation is simply the Spin$^c(n)$ representation $\Delta_n$.
\end{itemize}

We will also need the $Spin^{c,r}(n)$-equivariant map
\begin{eqnarray*}
 \mu_r\otimes\mu_n:\left(\ext^*\mathbb{R}^r\otimes_\mathbb{R} \ext^*\mathbb{R}^n\right)
 \otimes_\mathbb{R} (\Delta_r\otimes \Delta_n) &\longrightarrow& \Delta_r\otimes\Delta_n\\
(w_1 \otimes w_2)\otimes (\psi\otimes \varphi) &\mapsto& 
(w_1\otimes w_2)\cdot (\psi\otimes \varphi)
= (w_1\cdot\psi) \otimes (w_2\cdot \varphi).
\end{eqnarray*}

\subsection{Skew-symmetric 2-forms and endomorphisms associated to twisted spinors}

We will often write $f_{kl}$ for the Clifford product $f_k f_l$.

\begin{defi}
{\rm \cite{Espinosa-Herrera}}
Let $r\geq 2$, $\phi\in\Delta_r\otimes\Delta_n$, $X,Y\in\mathbb{R}^n$, $(f_1\ldots,f_r)$
an orthonormal basis of $\mathbb{R}^r$ and $1\leq k,l\leq r$.
\begin{itemize}
\item 
Define the real $2$-forms associated to the spinor $\phi$ by
\[\eta_{kl}^{\phi} (X,Y) = {\rm Re}\left< X\wedge Y\cdot \kappa_{r*}(f_kf_l)\cdot \phi,\phi\right>.\]

\item Define the antisymmetric endomorphisms
$\hat\eta_{kl}^\phi\in\End^-(\mathbb{R}^n)$ by
\[X\mapsto \hat\eta_{kl}^\phi(X):=(X\lrcorner \,\eta_{kl}^{\phi})^\sharp,\]
where $X\in\mathbb{R}^n$, $\lrcorner$ denotes contraction and $^\sharp$ denotes metric dualization from
$1$-forms to vectors.
\end{itemize}
\end{defi}

Observe that \(\eta_{kl}^\phi =  (\delta_{kl}-1)\eta_{lk}^\phi\). If $k\not= l$ then the imaginary part of $\eta_{kl}^{\phi}$ vanishes, so we can write
\[\eta_{kl}^{\phi} (X,Y) =\left< X\wedge Y\cdot \kappa_{r*}(f_kf_l)\cdot
\phi,\phi\right>.\]

\begin{lemma}{\rm \cite{Espinosa-Herrera}}
Any spinor $\phi\in\Delta_r\otimes\Delta_n$, $r\geq 2$, defines two maps (extended by linearity)
\begin{eqnarray*}
\ext^2 \mathbb{R}^r&\lra& \ext^2 \mathbb{R}^n\\
f_{kl} &\mapsto& \eta_{kl}^{\phi}
\end{eqnarray*}
and
\begin{eqnarray*}
\ext^2 \mathbb{R}^r&\lra& \End(\mathbb{R}^n)\\
f_{kl} &\mapsto& \hat\eta_{kl}^{\phi}.
\end{eqnarray*}
\end{lemma}

\subsection{Twisted partially pure spinors}

In order to simplify the statements, we will
consider the twisted Spin representation
\[
 \Sigma_r\otimes\Delta_n \subseteq \Delta_r\otimes\Delta_n.
\]
where 
\[
\Sigma_r=
\left\{
\begin{array}{ll}
\Delta_r & \mbox{if $r$ is odd,}\\
\Delta_r^+ & \mbox{if $r$ is even,} 
\end{array}
\right. 
\]
and $n,r\in\mathbb{N}$.

\begin{defi}
{\rm \cite{Herrera-Tellez}}
Let $(f_1,\ldots,f_r)$ be an  orthonormal frame of $\mathbb{R}^r$.
A unit-length spinor $\phi\in\Sigma_r\otimes\Delta_n$, $r<n$, is called a {\em twisted partially pure
spinor} if 
\begin{itemize}
 \item there exists a $(n-r)$-dimensional subspace $V^\phi\subset\mathbb{R}^n$ such that
for every $X\in V^\phi$, there exists a $Y\in V^\phi$ such that
\[X\cdot \phi = i\,\,Y\cdot \phi. \]

 \item it satisfies the equations
\begin{eqnarray*}
(\eta_{kl}^\phi + \kappa_{r*}(f_kf_l))\cdot \phi&=&0,\\
 \left<\kappa_{r*}(f_kf_l)\cdot \phi,\phi\right>&=&0,
\end{eqnarray*}
for all $1\leq k<l\leq r$.
\item If $r=4$, it also satisfies the condition
\[\left<\kappa(f_1f_2f_3f_4)\cdot \phi,\phi\right>=0.\]

\end{itemize}
\end{defi}

Let $(e_1,\ldots,e_{2m},e_{2m+1},\ldots,e_{2m+r})$ and 
$(f_1,\ldots,f_r)$ be orthonormal frames 
of $\mathbb{R}^{2m+r}$ and $\mathbb{R}^r$ res\-pec\-ti\-vely.
Consider the decomposition 
\[\Delta_{2m+r}=\Delta_r\otimes\Delta_{2m}^+\,\,\oplus\,\,\Delta_r\otimes\Delta_{2m}^-,\]
corresponding to the decomposition
\[\mathbb{R}^{2m+r}={\rm span}\{e_1,\ldots,e_{2m}\}
\oplus {\rm span}\{e_{2m+1},\ldots,e_{2m+r}\}.\]
Let
\[\varphi_0= u_{1,\ldots,1}\in \Delta_{2m}^+,\]
and
\[\{v_{\varepsilon_1,\ldots,\varepsilon_{[r/2]}}\,|\, (\varepsilon_1,\ldots,\varepsilon_{[r/2]}) \in
\{\pm1\}^{[r/2]}\}\]
be the unitary basis of the twisting factor $\Delta_r=\Delta({\rm span}(f_1,\ldots,f_r))$ which contains
$\Sigma_r$.
Let us define the standard twisted partially pure spinor
$\phi_0\in\Sigma_r\otimes\Delta_r\otimes\Delta_{2m}^+$ by
\begin{equation}
\phi_0 = \left\{
\begin{array}{ll}
{1\over \sqrt{2^{[r/2]}}}\,\, \left(\sum_{I\in\{\pm1\}^{\times [r/2]}}
v_I\otimes
\gamma_{r}(u_I)\right)\otimes \varphi_0 & \mbox{if $r$ is
odd,}\\
{1\over \sqrt{2^{[r/2]-1}}}\,\, \left(\sum_{I\in\left[\{\pm1\}^{\times [r/2]}\right]_+}
v_I\otimes
\gamma_{r}(u_I)\right)\otimes \varphi_0 &\mbox{if $r$ is even}, 
\end{array} \right.  \label{eq: canonical partilly pure spinor}
\end{equation}
where the elements of $\left[\{\pm1\}^{\times [r/2]}\right]_+$ 
contain an even number of $(-1)$.

We collect properties of partially pure spinors \cite{Herrera-Tellez} in the following proposition.

\vspace{.2in}

\begin{prop}
\begin{itemize} Let $\phi\in\Sigma_r\otimes \Delta_n$ be a partially pure spinor.
 \item 
 The definition of partially pure spinor does not depend on the choice of orthonormal basis of
$\mathbb{R}^r$.

\item
There exists an orthogonal complex structure on $V^\phi$ and $n-r \equiv 0$ {\rm (mod 2)} .

\item
If $r\geq2$,
\[\span\{\hat{\eta}_{kl}^\phi \in{\rm End}^-(\mathbb{R}^n) |\, 1\leq k< l\leq r\}\cong \mathfrak{so}(r).\]

\end{itemize}
 
\end{prop}
\qd

\subsection{Certain homogeneous spaces}\label{subsec: certain homogeneous spaces}

In this subsection, we present certain homogeneous spaces which will provide examples for 
various results in the following sections.

Consider the partial flag manifold 
\[\mathcal{G}_{m,s,r}={SO(2m+s+r)\over U(m)\times SO(s)\times SO(r)}\]
We will decompose the Lie algebra $\mathfrak{so}(2m+s+r)$
according to the natural inclusions
\[U(m)\times SO(s)\times SO(r)\subset SO(2m)\times SO(s)\times SO(r)\subset  SO(2m+s+r).\]
Note that
\begin{eqnarray*}
 \mathfrak{so}(2m+s+r) 
    &=& \ext^2 \mathbb{R}^{2m+s+r}\\
   &=& \ext^2 (\mathbb{R}^{2m}\oplus \mathbb{R}^s\oplus \mathbb{R}^r)\\
   &=& \ext^2 \mathbb{R}^{2m}\oplus \ext^2 \mathbb{R}^s\oplus \ext^2 \mathbb{R}^r 
   \oplus \mathbb{R}^{2m}\otimes \mathbb{R}^s\oplus \mathbb{R}^{2m}\otimes \mathbb{R}^r\oplus \mathbb{R}^s\otimes \mathbb{R}^s\\
   &=& \mathfrak{so}(2m) \oplus\mathfrak{so}(s) \oplus\mathfrak{so}(r) 
   \oplus \mathbb{R}^{2m}\otimes \mathbb{R}^s\oplus \mathbb{R}^{2m}\otimes \mathbb{R}^r\oplus \mathbb{R}^s\otimes \mathbb{R}^s\\
 \mathfrak{so}(2m)\otimes \mathbb{C} 
   &=&\ext^2(\mathbb{C}^m\oplus\overline{\mathbb{C}^m})\\
   &=&\ext^2\mathbb{C}^m\oplus\mathbb{C}^m\otimes\overline{\mathbb{C}^m}\oplus\ext^2\overline{\mathbb{C}^m}\\
   &=&[[\ext^2\mathbb{C}^m]]\otimes \mathbb{C}\oplus\mathfrak{u}(m)\otimes
    \mathbb{C},\\
 \mathbb{R}^{2m}&=& [[\mathbb{C}^m]],
\end{eqnarray*}
where the symbol $[[\mathbb{C}^m]]$ denotes the underlying real vector space $\mathbb{R}^{2m}$ of $\mathbb{C}^m$ carrying a complex structure.
Thus
\[\mathfrak{so}(2m+s+r) = \mathfrak{u}(m) \oplus \mathfrak{so}(s) \oplus \mathfrak{so}(r) \oplus 
\left( [[\ext^2\mathbb{C}^m]] 
\oplus [[\mathbb{C}^m]]\otimes \mathbb{R}^s \oplus [[\mathbb{C}^m]]\otimes \mathbb{R}^r 
\oplus \mathbb{R}^s\otimes \mathbb{R}^r 
\right)
\]
and the tangent space of $\mathcal{G}_{m,s,r}$ decomposes as follows
\[T_{{\rm Id}}\mathcal{G}_{m,s,r} \cong  
 [[\ext^2\mathbb{C}^m]] \oplus [[\mathbb{C}^m]]\otimes \mathbb{R}^s\oplus [[\mathbb{C}^m]]\otimes \mathbb{R}^r \oplus\mathbb{R}^s\otimes \mathbb{R}^r. 
\]
This gives the isotropy representation
\begin{eqnarray*}
 U(m)\times SO(s)\times SO(r) &\longrightarrow& SO(T_{{\rm Id}}\mathcal{G}_{m,s,r})\\
 (A,B,C) &\mapsto& \left(
 \begin{array}{llll}
[[\ext^2A]] &  &  & \\
 & [[A]]\otimes B &  & \\
 &  & [[A]]\otimes C & \\
 &  &  & B\otimes C
 \end{array}
\right),
\end{eqnarray*}
where $\ext^2A$ denotes the linear transformation induced by $A$ on $\ext^2\mathbb{C}^m$, $[[A]]$ the transformation $A$ viewed as a real linear transformation on 
$[[\mathbb{C}^m]]=\mathbb{R}^{2m}$, and $B\otimes C$ the induced transformation on $\mathbb{R}^s\otimes\mathbb{R}^r$ (i.e. the Kronecker product of $B$ and $C$).

\section{Doubly twisted Spin structures}\label{sec: doubly twisted Spin structures}

In this section, we introduce the (doubly) twisted Spin structures we need to carry out our spinorial characterization of CR structures, 
and the corresponding twisted Dirac operator and Laplacian. We deduce some topological conditions on manifolds that support such  structures,  
a Schr\"odinger-Lichnerowicz type formula, and give some Bochner-type arguments.

\begin{defi}\label{defi-twisted-spin-structure}
Let $M$ be an oriented $n$-dimensional Riemannian manifold, $P_{SO(M)}$ be its principal bundle of
orthonormal frames and $r\in\mathbb{N}$.
% , $r\geq 2$. 
A Spin$^{c,r}(n)$ structure on $M$ consists of 
an auxiliary principal $SO(r)$ bundle $P_{SO(r)}$, 
an auxiliary principal $U(1)$ bundle $P_{U(1)}$ 
and
a principal Spin$^{c,r}(n)$ bundle $P_{Spin^{c,r}(n)}$ together with an equivariant $2:1$ covering map
\[\Lambda:P_{Spin^{c,r}(n)}\lra P_{SO(M)} \tilde{\times} P_{SO(r)}\tilde{\times} P_{U(1)},\]
where $\tilde{\times}$ denotes the fibered product,
such that $\Lambda(pg)=\Lambda(p)(\lambda_{n,r,2})(g)$ for all $p\in
P_{Spin^{c,r}(n)}$ and
$g\in Spin^{c,r}(n)$, where $\lambda_{n,r,2}:Spin^{c,r}(n)\lra SO(n)\times SO(r)\times
U(1)$ denotes the
canonical $2$-fold cover. 

A $n$-dimensional Riemannian manifold $M$ admitting a Spin$^{c,r}(n)$ structure will be called a {\em
Spin$^{c,r}$ manifold}.
\end{defi}

{\bf Remark}. 
A Spin$^{c,r}$ manifold with trivial $P_{SO(r)}$ and $P_{U(1)}$ auxiliary bundles is a Spin manifold.
On the other hand, we have the following:
\begin{itemize}
 \item 
Any Spin manifold admits a Spin$^{c,r}$ structure with trivial $P_{SO(r)}$ and $P_{U(1)}$ auxiliary
bundles via the inclusion $Spin(n)\subset Spin^{c,r}(n)$.

 \item 
Any Spin$^c$ manifold admits a Spin$^{c,r}$ structure with trivial $P_{SO(r)}$ auxiliary
bundle via the inclusion $Spin^c(n)\subset Spin^{c,r}(n)$.

 \item 
Any Spin$^r$ manifold (cf. \cite{Espinosa-Herrera}) admits a Spin$^{c,r}$ structure with trivial $P_{U(1)}$ auxiliary
bundle via the inclusion $Spin^r(n)\subset Spin^{c,r}(n)$.

\end{itemize}

\subsection{Existence of Spin$^{c,r}$ structures}\label{subsec: existence of non-reducible twisted Spin structures}

We will characterize the existence of a Spin$^{c,r}$ structure in terms of a Spin structure. 

\begin{prop}{\rm \cite[p. 47]{Friedrich}}
Let $G\subset SO(N)$ be a connected compact Lie subgroup with 
\[\pi_1(SO(N)/G)=\{0\}.\]
A $G$-principal bundle $Q$ over a connected CW-complex $X$ has a Spin structure if and only if there exists a homomorphism 
$f:\pi_1(Q)\to\pi_1(SO(N))$ for which the diagram
\[
\xymatrix{
\pi_1(G)\ar[d]_h \ar[r]^-{i_{\#}} & \pi_1(SO(n))\\
\pi_1(Q) \ar[ur]_{f}
}
\]
commutes.
\qd
\end{prop}

By setting $N=n+r+2$, $G=SO(n)\times SO(r)\times U(1)$, $Q=P_{SO(M)} \tilde{\times} P_{SO(r)}\tilde{\times} P_{U(1)}$ and 
considering the natural inclusion of $SO(n)\times SO(r)\times U(1)\subset SO(n+r+2)$ 
we have that 
\[\pi_1\left({SO(n+r+2)\over SO(n)\times SO(r)\times U(1)}\right)=\{0\}.\]

\begin{corol}\label{pgsub}
The bundle $P_{SO(M)} \tilde{\times} P_{SO(r)}\tilde{\times} P_{U(1)}$
over $M$ has a Spin structure if and only if there exists a homomorphism $f:\pi_1(Q)\to\pi_1( SO(n+r+2))$ for which 
the diagram
\[
\xymatrix{
\pi_1(SO(n)\times SO(r)\times U(1))\ar[d]_h \ar[r]^-{i_{\#}} & \pi_1(SO(n+r+2))\\
\pi_1(P_{SO(M)} \tilde{\times} P_{SO(r)}\tilde{\times} P_{U(1)}) \ar[ur]_{f}
}
\]
commutes. 
\qd
\end{corol}

\begin{lemma} \label{lemma: calculo exactitud}For $r\geq 2$,
\[\ker(i_{\#})= (\lambda_n\times\lambda_r\times \lambda_2)_{\#}(\pi_1(Spin^{c,r}(n))).\] 
\end{lemma}
{\em Proof}. 
Using additive notation, we have for $r\geq 3$
\begin{eqnarray*}
 \pi_1(SO(n)\times SO(r)\times U(1))\cong \mathbb{Z}_2\oplus\mathbb{Z}_2\oplus\mathbb{Z} &\xrightarrow{i_{\#}}& \pi_1(SO(n+r+2))\cong \mathbb{Z}_2\\
(a,b,c) &\mapsto& a+b+c \,\,\,\mbox{\rm (mod 2)},
\end{eqnarray*}
and for $r=2$
\begin{eqnarray*}
 \pi_1(SO(n)\times SO(2)\times U(1))\cong \mathbb{Z}_2\oplus\mathbb{Z}\oplus\mathbb{Z} &\xrightarrow{i_{\#}}& \pi_1(SO(n+4))\cong \mathbb{Z}_2\\
(a,b,c) &\mapsto& a+b+c \,\,\,\mbox{\rm (mod 2)}.
\end{eqnarray*}

For $r\geq 3$, we see that $\pi_1(Spin^{c,r}(n))\cong \mathbb{Z}_2\oplus\mathbb{Z}$, as described in Lemma \ref{lemma: fundamental group SpinCR}, is mapped as follows
\[\begin{array}{ccccc}
 \pi_1(Spin^{c,r}(n)) &\xrightarrow{(\lambda_n\times\lambda_r\times \lambda_2)_{\#}}& 
 \pi_1(SO(n)\times SO(r)\times U(1)) &\xrightarrow{i_{\#}}&\pi_1(SO(n+r+2)\\
 (a,b)&\mapsto& (a,a+b ,b) &\mapsto& 2a+2b=0 \,\,\mbox{\rm (mod 2)}.
\end{array}
\]

For $r=2$, we see that $\pi_1(Spin^{c,2}(n))\cong \mathbb{Z}\oplus\mathbb{Z}$, as described in Lemma \ref{lemma: fundamental group SpinCR}, is mapped as follows
\[\begin{array}{ccccc}
 \pi_1(Spin^{c,2}(n)) &\xrightarrow{(\lambda_n\times\lambda_2\times \lambda_2)_{\#}}& 
 \pi_1(SO(n)\times SO(2)\times U(1)) &\xrightarrow{i_{\#}}&\pi_1(SO(n+4)\\
 (a,b)&\mapsto& (a+b,a ,b) &\mapsto& 2a+2b=0 \,\,\mbox{\rm (mod 2)}.
\end{array}
\]
\qd

\begin{prop}\label{equivfb}
$M$ admits an $SO(r)\times SO(2)$-principal bundle $P_{SO(r)\times SO(2)}$ such that the fibre 
product $Q=P_{SO(n)}\tilde{\times}P_{SO(r)\times SO(2)}$ has a Spin structure if and only if $M$ has a Spin$^{c,r}$ structure.
\end{prop}

{\em Proof}.
If $M$ has a Spin$^{c,r}$ structure then $P_{SO(2)}:=P/Spin^r(n)$ and $P_{SO(r)}:=P/Spin^c(n)$ are $SO(2)$ and $SO(r)$ principal bundles respectively, 
so that  $P_{SO(r)\times SO(2)}:=P_{SO(r)}\tilde\times P_{SO(2)}$ is a 
$SO(r)\times SO(2)$ principal bundle over $M$. Now, there exists an injective homomorphism $\tilde{i}$ which makes the diagram
\[
\xymatrix{
Spin^{c,r}(n)\ar[r]^-{\tilde{i}}\ar[d] & Spin(n+r+2)\ar[d]\\
SO(n)\times SO(r)\times SO(2)\ar[r]_-{i}&  SO(n+r+2)\\
}
\]
commute. From this we obtain a Spin structure for $Q$ in the sense of Corollary $\ref{pgsub}$.

Conversely, 
let $\lambda=\lambda_n\times\lambda_r\times \lambda_2$ and $F=SO(n)\times SO(r)\times U(1)$. According to Corollary \ref{pgsub},
due to the existence of $f$, $H=\ker(f)\subset \pi_1(Q)$ is a subroup of index $2$. Therefore, 
there exists a double covering space $\Lambda:P_{Spin^{c,r}(n)}\to Q$ 
corresponding to $H$.
Let $\mu:Q\times F\to Q$ be the action of $F$ in $Q$ and consider the composition 
of induced maps on fundamental groups 
$$
\xymatrix{
 \pi_1(P_{Spin^{c,r}(n)}\times Spin^{c,r}(n))\ar[r]^-{(\Lambda\times\lambda)_{\#}} & \pi_1(Q\times F)\ar[r]^-{\mu_{\#}} & \pi_1(Q). 
}
$$
If $(\sigma,\tau)\in\pi_1(P_{Spin^{c,r}(n)})\times\pi_1(Spin^{c,r}(n))$, by means of the inclusion $h$,
\begin{eqnarray*}
\mu_{\#}\circ(\Lambda\times\lambda)_{\#}(\sigma,\tau)
   &=&
  \Lambda_{\#}(\sigma)\lambda_{\#}(\tau)\\ 
   &=&
  \Lambda_{\#}(\sigma)*h(\lambda_{\#}(\tau)) 
\end{eqnarray*}
where $*$ denotes product in the relevant fundamental group.
We know that 
\[\Lambda_{\#}(\sigma)\in H \quad\quad\mbox{and} \quad\quad f(h(\lambda_{\#}(\tau)))=i_{\#}(\lambda_{\#}(\tau))=0\]
by Lemma \ref{lemma: calculo exactitud} and Corollary \ref{pgsub}. Thus,  $h(\lambda_{\#}(\tau))\in H$ and  
$\Lambda_{\#}(\sigma)*h(\lambda_{\#}(\tau))\in H$.
Hence, there exists a lift $\tilde{\mu}:P_{Spin^{c,r}(n)}\times Spin^{c,r}(n)\to P_{Spin^{c,r}(n)}$ 
which gives the equivariance in Definition \ref{defi-twisted-spin-structure}. 
\qd

Now, we will derive a condition for a 
simply connected manifold to have a ``non-reducible'' Spin$^{c,r}$ structure, i.e. a Spin$^{c,r}$ structure 
which does not come from a Spin, nor a Spin$^{c}$, nor Spin$^r$ structure.\\

\begin{prop}\label{propscr}
Let $M$ be simply connected and $Q$ its $SO(n)$-principal bundle of orthonormal frames. 
The following are equivalent
\begin{enumerate}
 \item $Q$ has a Spin$^{c,r}$ structure but does not have a Spin, nor a Spin$^c$, nor a Spin$^r$ structure.
 \item There exists a $SO(r)\times SO(2)$ bundle $P_1$ over $X$ such that in the long exact sequence
 $$
 \xymatrix{
 \cdots\ar[r] & \pi_2(X)\ar[r]^-{\partial} & \pi_1(SO(n)\times SO(r)\times  SO(2))\ar[r]^-h & \pi_1(Q\tilde{\times}P_1)\ar[r] & \pi_1(X)=0, 
 }
 $$
$\Ima(\partial)\cong\langle(1,0,p),(0,1,p)\rangle\subset\mathbb{Z}_2\oplus\mathbb{Z}_2\oplus\mathbb{Z}$ with $p$ odd.
\end{enumerate}
\end{prop}

{\em Proof}. 
If $(P,\Lambda)$ is a Spin$^{c,r}$ structure on $Q$ then $P_{SO(2)}:=P/Spin^r(n)$ and $P_{SO(r)}:=P/Spin^c(n)$ 
are $ SO(2)$ and $ SO(r)$ principal bundles respectively, so that  $P_1:=P_{ SO(r)}\tilde\times P_{ SO(2)}$ is a 
$SO(r)\times SO(2)$ principal bundle over $X$. Now, by Proposition $\ref{equivfb}$, the fibre product $Q\tilde{\times}P_1$ 
has a Spin structure. 
By Corollary $\ref{pgsub}$, this means that there exists a map 
$f:\pi_1(Q\tilde{\times}P_1)\to\pi_1(SO(n+r+2))$ such that the diagram
$$
\xymatrix{
\pi_1(SO(n)\times SO(r)\times SO(2))=\mathbb{Z}_2\oplus\mathbb{Z}_2\oplus\mathbb{Z}\ar[d]_{h} \ar[r]^-{i_{\#}} & \pi_1(SO(n+r+2))=\mathbb{Z}_2\\
\pi_1(Q\tilde{\times}P_1) \ar[ur]_{f}
}
$$
commutes.
Now, if $Q$ does not have a Spin structure then we have $\pi_1(Q)=0$ in the following commutative diagram 
% \begin{gather}
% \begin{aligned}
% \xymatrix{
% & \vdots\ar[d]\\
% & \pi_2(X) \ar[d]^{\partial}\\
% & \mathbb{Z}_2\oplus\mathbb{Z}_2\oplus\mathbb{Z}\ar[d]_h\\
% \mathbb{Z}_2\oplus\mathbb{Z}=\pi_1(SO(r)\times SO(2)) \ar[r]^-{k}\ar[ur]^{j_{\#}}& \pi_1(Q\tilde{\times}P_1)\ar[r] \ar[d] & \pi_1(Q)\ar[r] & 0\\
% &\pi_1(X)=0.
% }
% \end{aligned}
% \label{diagram1}
% \end{gather}
\begin{equation}
\xymatrix{
& \vdots\ar[d]\\
& \pi_2(X) \ar[d]^{\partial}\\
& \mathbb{Z}_2\oplus\mathbb{Z}_2\oplus\mathbb{Z}\ar[d]_h\\
\mathbb{Z}_2\oplus\mathbb{Z}=\pi_1(SO(r)\times SO(2)) \ar[r]^-{k}\ar[ur]^{j_{\#}}& \pi_1(Q\tilde{\times}P_1)\ar[r] \ar[d] & \pi_1(Q)\ar[r] & 0\\
&\pi_1(X)=0.
}
\label{diagram1}
\end{equation}
Thus, $k$ is onto and 
$$(\mathbb{Z}_2\oplus\mathbb{Z}_2\oplus\mathbb{Z})/\Ima(\partial)\cong\pi_1(Q\tilde{\times}P_1)
=h(\mathbb{Z}_2\oplus\mathbb{Z}_2\oplus\mathbb{Z})=k(\mathbb{Z}_2\oplus\mathbb{Z}).$$
Now, we will describe the group $K=\pi_1(Q\tilde{\times}P_1)$. 
It depends on the nontrivial elements $h(1,0,0)=\alpha$, $h(0,1,0)=\beta$ and $h(0,0,1)=\gamma$. 
First, 
we have $K=\langle \beta, \gamma\rangle$, so that
$\alpha=a\beta+b\gamma$ for some integers $a,\,b$. Since $\alpha$ and $\beta$ have order two in $K$
\[0=2\alpha=2a\beta+2b\gamma=2b\gamma,\]
and $K$ is a finite group.
Now,
\begin{itemize}
 \item[(i)] If $\beta\in\langle \gamma \rangle$ then $\gamma$ has order $2p$, $K\cong\mathbb{Z}_{2p}$ and 
$\alpha=\beta=p\gamma$. Since there is only one nontrivial map $f:\mathbb{Z}_{2p}\to\mathbb{Z}_2$, 
$f\circ h=i_{\#}$ if and only if $p$ is odd. 

 \item[(ii)] If $\beta\notin\langle \gamma \rangle$ then $K\cong\mathbb{Z}_2\oplus\mathbb{Z}_{d}$. 
If $d$ is odd, $f:\mathbb{Z}_d\to\mathbb{Z}_2$ must be trivial and $(f\circ h)(0,0,1)=0$ which 
gives us no Spin$^{c,r}$ structure. If $d=2p$, then
\begin{equation*}
\alpha= \beta, \quad\mbox{or}\quad \alpha=p\gamma \quad\mbox{or}\quad \alpha=\beta+p\gamma. \label{eq: 3 options}
\end{equation*}
In order to have $i_{\#}=f\circ h$ and, therefore, the existence of the Spin$^{c,r}$ structure, if $\alpha = p\gamma$
then $p$ must be odd, and if $\alpha=\beta+p\gamma$ then $p$ must be even.

\end{itemize}
Now, we are going to rule out the three options in (ii). Note that $K'=\pi_1(Q\tilde\times P_{SO(r)}\tilde\times P_{SO(2)})$
and the $SO(2)$ fibre bundle $Q\tilde\times P_{SO(r)}\tilde\times P_{SO(2)}\to Q\tilde\times P_{SO(r)}$ 
gives the commutative diagram
\begin{equation}
\xymatrix{
& & \vdots\ar[d] & \vdots\ar[d]\\
& & \pi_2(X) \ar[d] &  \pi_2(X) \ar[d]\\
\cdots\ar[r] & \mathbb{Z}=\pi_1(SO(2))\ar[r]& \mathbb{Z}_2\oplus\mathbb{Z}_2\oplus\mathbb{Z}\ar[r]\ar[d] & \mathbb{Z}_2\oplus\mathbb{Z}_2\ar[d]\ar[r] & 0\\
\cdots\ar[r] &\mathbb{Z}=\pi_1(SO(2)) \ar[r]\ar[ur]& \pi_1(Q\tilde{\times}P_{SO(r)}\tilde\times P_{SO(2)})\ar[r] \ar[d] & \pi_1(Q\tilde\times P_{SO(r)})\ar[r]\ar[d] & 0\\
& &\pi_1(X)=0 & \pi_1(X)=0
}
\label{diagram2}
\end{equation}
\begin{itemize}
 \item If $K'=\mathbb{Z}_2\oplus\mathbb{Z}_{2p}=\langle\beta,\,\gamma\rangle$ with $\alpha=\beta$ then, 
by exactness of the diagram, $\pi_1(Q\tilde\times P_{SO(r)})=(\mathbb{Z}_2\oplus\mathbb{Z}_{2p})/\langle\gamma\rangle=\mathbb{Z}_2
\cong \langle \beta \rangle$, which gives us a Spin$^r$ structure. 

 \item The same happens if $\alpha=\beta+p\gamma$. 
The quotient is isomorphic to $\mathbb{Z}_2$, whose equivalence classes are
\[\{(0,0),(0,1),\ldots,(0,2p-1)\}
\quad\mbox{and}\quad
\{(1,0),\ldots,(1,2p-1)\},\]
where $\gamma=(0,1)$ belongs to the first one, and $\alpha=(1,p)$ and $\beta=(1,0)$ belong to the second one.
In other words,
$\alpha$ and $\beta$ are mapped to the nontrivial class and we have a Spin$^r$ structure.

 \item Now if $K'=\mathbb{Z}_2\oplus\mathbb{Z}_{2p}$ with $p$ odd and $\alpha=p\gamma$, the $SO(r)$ fibre 
bundle $Q\tilde\times P_{SO(r)}\tilde\times P_{SO(2)}\to Q\tilde\times P_{SO(2)}$ gives the commutative diagram
$$
\xymatrix{
& & \vdots\ar[d] & \vdots\ar[d]\\
& & \pi_2(X) \ar[d] &  \pi_2(X) \ar[d]\\
\cdots\ar[r] & \mathbb{Z}_2=\pi_1(SO(r))\ar[r]& \mathbb{Z}_2\oplus\mathbb{Z}_2\oplus\mathbb{Z}\ar[r]\ar[d] & \mathbb{Z}_2\oplus\mathbb{Z}\ar[d]\ar[r] & 0\\
\cdots\ar[r] &\mathbb{Z}_2=\pi_1(SO(r)) \ar[r]\ar[ur]& \pi_1(Q\tilde{\times}P_{SO(r)}\tilde\times P_{SO(2)})\ar[r] \ar[d] & \pi_1(Q\tilde\times P_{SO(2)})\ar[r]\ar[d] & 0\\
& &\pi_1(X)=0 & \pi_1(X)=0
}
$$
so that 
$\pi_1(Q\tilde\times P_{SO(2)})=\mathbb{Z}_2\oplus\mathbb{Z}_{2p}/\langle\beta\rangle=\mathbb{Z}_{2p}$ 
with $p$ odd, which implies the existence of a Spin$^c$ structure.

\end{itemize}

Now we know that $K=\mathbb{Z}_{2p}=\langle\gamma\rangle$ with $p$ odd and $\alpha=\beta=p\gamma$. 
This should be the same as the quotient $(\mathbb{Z}_2\oplus\mathbb{Z}_2\oplus\mathbb{Z})/\Ima(\partial)$ where, 
by extacness, $\Ima(\partial)=\ker(h)$. We see that the map $h$ is given by $h(a,b,c)=((a+b)p+c)\gamma$. 
The kernel of this map is given by the $(a,b,c)\in\mathbb{Z}_2\oplus \mathbb{Z}_2\oplus \mathbb{Z}$ such 
that $(a+b)p+c\equiv 0$ (mod $2p$), i.e. $(a,b,c)\in\langle (0,1,p), (1,0,p)\rangle$ where $p$ is odd.

Conversely, assume $\Ima(\partial)=\langle(1,0,p),(0,1,p)\rangle$, $p$ odd, and put this in the diagram (\ref{diagram1}). 
By exactness of the column, $\pi_1(Q\tilde{\times}P_1)\cong\mathbb{Z}_{2p}$. This group is generated by 
the non trivial element $\gamma=h(0,0,1)$, and we have $h(1,0,0)=h(0,1,0)=p\gamma$. Thus,
$k=h\circ j_{\#}$ is onto, we have no Spin structure and the only nonzero homomorphism $f:\mathbb{Z}_{2p}\to\mathbb{Z}_2$ 
gives us $i_{\#}=f\circ h$, i.e. the existence of a  Spin$^{c,r}$ structure.

The $SO(r)$ bundle $P_{SO(r)}=P_1/SO(2)$ fits into a commutative diagram similar to (\ref{diagram2}). 
By exactness, we have $\pi_1(Q\tilde\times P_{SO(r)})=\{0\}$ and we cannot have a Spin$^r$ structure. 
Similarly, the $SO(2)$ bundle $P_{SO(2)}=P_1/SO(r)$ fits into a similar diagram
so that $\pi_1(Q\tilde\times P_{SO(2)})=\mathbb{Z}_p$, 
and there is not map $f$ as in Corollary \ref{pgsub} to have a Spin$^c$ structure.
Note that in the case $p=1$ this last group is zero.
\qd

{\bf Example.} Now we will give an example of a manifold satisfying the conditions of the previous Proposition. 
Let $X=G/H$ with $G=SO(2m+2+r)$ and $H=U(m)\times U(1)\times SO(r)$, $r\geq 3$. Since $H$ is a compact connected
subgroup of $G$ and the inclusion map induces a map of fundamental groups which is onto,  $\pi_1(X)=\{0\}$.

Now consider the bundle of orthonormal frames $Q=G\times_{\rho}SO(n)$ where $n=m^2+2mr+3m+2r=\dim(X)$ 
and 
\[\rho: H\hookrightarrow G\times SO(n)\]
is given by the inclusion of $H$ into the first factor and the isotropy 
representation in the second
which is given by 
\begin{eqnarray*}
\xi:H&\longrightarrow& SO(n)\\
(A,e^{i\theta},B)&\mapsto& [[\Lambda^2 A]]\oplus([[A]]\otimes R_{\theta})\oplus([[A]]\otimes B)\oplus(R_{\theta}\otimes B), 
\end{eqnarray*}
where $R_{\theta}$ is the rotation in $\mathbb{R}^2$ by an angle of $\theta$.
This gives a fibration 
\[\begin{array}{ccc}
H & \hookrightarrow & G\times SO(n)\\
 &  & \downarrow\\
 &  & Q
  \end{array}
\]
which induces the long exact sequence of homotopy groups
\[
\xymatrix{
\cdots\ar[r] & \pi_1(H)=\mathbb{Z}\oplus\mathbb{Z}\oplus\mathbb{Z}_2\ar[r]^-{\rho_{\#}}& \pi_1(G\times SO(n))=\mathbb{Z}_2\oplus\mathbb{Z}_2\ar[r]& \pi_1(Q)\ar[r] & 0.\\
}
\]
First note that the isotropy representation induces the map 
\begin{eqnarray*}
\pi_1(H) &\longrightarrow& \pi_1(SO(n))\\
(a,b,c)&\mapsto&(m-1+2+r)a+(2m+r)b+(2m+2)c\,\,\,\,\,\mbox{(mod 2)}\\ 
   &=&(m-1+r)a+rb\,\,\,\,\,\mbox{(mod 2)}.
\end{eqnarray*}
Thus,
\[\rho_{\#}(a,b,c)=((a+b+c)\,\,\mbox{(mod 2)}, ((m-1+r)a+br)\,\,\mbox{(mod 2)}).\]
Note that $(m-1+r)a+br\equiv 0$ (mod 2) if and only if $r$ is even and $m$ is odd. So, by exactness, 
$\pi_1(Q)=\mathbb{Z}_2$ ($Q$ has a Spin structure) if and only if $r$ is even and $m$ is odd.

Let $m$ and $r$ be even and consider 
\begin{eqnarray*}
\sigma:H&\longrightarrow& (G\times SO(n))\times U(1)\times SO(r)\\
(A,e^{i\theta},B)&\mapsto& (\rho(A,e^{i\theta},B),e^{i\theta},B)
\end{eqnarray*}
We have the fibration 
$H\hookrightarrow G\times SO(n)\times U(1)\times SO(r)\xrightarrow{\nu} G\times_{\sigma} SO(n)\times U(1)\times SO(r)=Q\tilde{\times} P_{U(1)}\tilde{\times}P_{SO(r)}$ 
which gives
\[
\xymatrix{
\cdots\ar[r] & \mathbb{Z}\oplus\mathbb{Z}\oplus\mathbb{Z}_2\ar[r]^-{\sigma_{\#}}
& \mathbb{Z}_2\oplus\mathbb{Z}_2\oplus\mathbb{Z}\oplus\mathbb{Z}_2\ar[r]^{\nu_{\#}}
& \pi_1(Q\tilde{\times} P_{U(1)}\tilde{\times}P_{SO(r)})\ar[r] & 0,\\
}
\]
where 
\[\sigma_{\#}(a,b,c)=((a+b+c)\,\,\mbox{(mod 2)}, a\,\,\mbox{(mod 2)},b,c\,\,\mbox{(mod 2)}).\]
We see that $\Ima(\sigma_{\#})=\langle (1,1,0,0),\,(1,0,1,0),\,(1,0,0,1)\rangle=L$ is a subgroup of index two 
because $(1,0,0,0)\notin L$ 
and $((1,0,0,0)+L)\cup L=\mathbb{Z}_2\oplus\mathbb{Z}_2\oplus\mathbb{Z}\oplus\mathbb{Z}_2$. 
By exactness, $\pi_1(Q\tilde{\times} P_{U(1)}\tilde{\times}P_{SO(r)})\cong\mathbb{Z}_2$.

Consider $f:\pi_1(Q\tilde{\times} P_{U(1)}\tilde{\times}P_{SO(r)})\to\pi_1(SO(n+2+r))$ to be the only nontrivial 
homomorphism between these groups. Now, the inclusion of the fiber $SO(n)\times Un(1)\times SO(r)$ into the fiber bundle 
$Q\tilde{\times} P_{U(1)}\tilde{\times}P_{SO(r)}$ is given by 
the inclusion $j$ into the last three factors of 
$G\times SO(n)\times U(1)\times SO(r)$, followed by the projection $\nu$. Thus, the map $h$ in Proposition $\ref{propscr}$ 
is given by $h=\nu_{\#}\circ j_{\#}$.

Consider, for simplicity, $\pi_1(Q\tilde{\times} P_{U(1)}\tilde{\times}P_{SO(r)})=\{0,1\}$. 
From the explicit description of $L$, we can see 
\[h(a,b,c)=\left\{\begin{array}{ll}
0 & \mbox{if $a+b+c\equiv 0$ (mod 2)},\\
1 & \mbox{if $a+b+c\equiv 1$ (mod 2)}.
                  \end{array}
\right.\]
This means that $f\circ h$  is the same map as the inclusion of 
$\pi_1(SO(n)\times U(1)\times SO(r))\hookrightarrow \pi_1(SO(n+2+r))$. 
By Proposition $\ref{propscr}$, $X$ has a Spin$^{c,r}$ structure which does not come from either a Spin, nor a 
Spin$^c$, nor a Spin$^r$ structure.
\qd

\subsection{Covariant derivatives and twisted differential operators}

Let $M$ be a Spin$^{c,r}$ $n$-dimensional manifold,
 $\omega$ the Levi-Civita connection 1-form on its principal bundle of orthonormal frames $P_{SO(n)}$,
 $\theta$ and $iA$ chosen  connection 1-forms on the auxiliary bundles $P_{SO(r)}$ and $P_{U(1)}$ respectively.
These connections forms give rise to covariant derivatives $\nabla$, $\nabla^\theta$ and $\nabla^A$ on the associated vector bundles 
\begin{eqnarray*}
TM &=& P_{Spin^{c,r}(n)}\times_{\lambda_{n,r,2}} (\mathbb{R}^n\times\{0\}\times\{0\}),\\
F&=&P_{Spin^{c,r}(n)}\times_{\lambda_{n,r,2}} (\{0\}\times\mathbb{R}^r)\times\{0\}),\\
L&=&P_{Spin^{c,r}(n)}\times_{\lambda_{n,r,2}} (\{0\}\times\{0\}\times\mathbb{C}),
\end{eqnarray*}
Furthermore, the three connections help define a connection on the twisted spinor bundle
\begin{eqnarray*}
S&=& P_{Spin^{c,r}(n)}\times_{\kappa_n^{c,r}} 
(\Sigma_r\otimes\Delta_n)
\end{eqnarray*}
given (locally) as follows
\begin{eqnarray*}
\nabla^{\theta,A}&:& \Gamma(S)\lra \Gamma(T^*M\otimes
S)\\
\nabla^{\theta,A} (\varphi\otimes \psi) 
   &=&
  d(\varphi\otimes\psi) +\varphi\otimes\left[{1\over 2}\sum_{1\leq i<j\leq n}\omega_{ij}\otimes
e_ie_j\cdot\psi\right] \\
   &&
+\left[{1\over 2}\sum_{1\leq k<l\leq r}\theta_{kl}\otimes \kappa_{r*}(f_kf_l)\cdot
\varphi\right]\otimes \psi
+{i\over 2} \,\varphi \otimes (A\cdot\psi),
\end{eqnarray*}
where $\varphi\otimes\psi\in\Gamma(S)$, 
$(e_1,\ldots,e_n)$ and $(f_1,\ldots,f_r)$ are a local orthonormal frames of $TM$ and $F$
resp., 
$\omega_{ij}$, $\theta_{kl}$ and $A$ are the corresponding local connection 1-forms for $TM$, $F$ and $L$ respectively.

From now on, we will omit the upper and lower bounds on the indices, by declaring $i$ and $j$
to be the indices for the frame vectors of $TM$, and $k$ and $l$ to be the indices for the frame sections of $F$. 

Now, for any tangent vectors $X,Y\in T_xM$, the spinorial curvature is defined by
\begin{eqnarray}
R^{\theta,A}(X,Y)(\varphi\otimes \psi)
   &=&
  \varphi\otimes\left[{1\over 2} \sum_{i<j} \Omega_{ij}(X,Y) e_ie_j\cdot\psi\right]\nonumber\\
   &&
  +\left[{1\over 2} \sum_{k<l} \Theta_{kl}(X,Y) \kappa_{r*}(f_kf_l)\cdot\varphi\right]\otimes\psi
  +{i\over 2}\,\varphi\otimes(dA(X,Y)\psi) 
,
\label{curvature-twisted-spin}
\end{eqnarray}
where
\[\Omega_{ij}(X,Y)=\big<R^M(X,Y)(e_i),e_j
\big>\quad\quad\mbox{and}\quad\quad\Theta_{kl}(X,Y)=\big<R^F(X,Y)(f_k),f_l \big>.\]
Here $R^M$ (resp. $R^F$) denotes the curvature tensor of $M$ (resp. of $F$).

For $X, Y$ vector fields and $\phi\in\Gamma(S)$ a spinor field,
we have compatibility of the covariant derivative with Clifford multiplication, 
\[
\nabla^{\theta,A}_X(Y\cdot\phi) = (\nabla_XY)\cdot\phi +
Y\cdot\nabla_X^{\theta,A}\phi.
\]

\begin{defi}
The {\em twisted Dirac operator} is the first order differential operator 
$\dirac^{\theta,A}:\Gamma(S)\longrightarrow \Gamma(S)$
defined by
\begin{eqnarray*}
\dirac^{\theta,A}(\phi)
&=&\sum_{i=1}^n e_i\cdot\nabla_{e_i}^{\theta,A}(\phi).
\end{eqnarray*}
\end{defi}

{\bf Remark}.
The twisted Dirac operator $\dirac^{\theta,A}$ is well-defined and formally self-adjoint on compact
manifolds. Moreover, 
if $h\in C^{\infty}(M)$, $\phi\in\Gamma(S)$, we have
\[
\dirac^{\theta,A}(h\,\phi) = {\rm grad}(h)\cdot \phi +
h\,\dirac^{\theta,A}(\phi).
\]
The proofs of these facts are analogous to the ones for the Spin$^c$ Dirac operator
\cite{Friedrich}.

\begin{defi}
The {\em twisted Spin connection Laplacian} is the second order differential operator 
$\Delta^{\theta,A}: \Gamma(S)\rightarrow\Gamma(S)$  defined as
\[
\Delta^{\theta,A}(\phi) = - \sum_{i = 1}^n
\nabla^{\theta,A}_{e_i}\nabla^{\theta,A}_{e_i}(\phi) - \sum_{i = 1}^n
{\rm div}(e_i)\nabla^{\theta,A}_{e_i}(\phi). 
\]
\end{defi}

\subsection{A Schr\"odinger-Lichnerowicz-type formula}\label{section-twisted-SL-formula}

Just as in \cite{Friedrich,Espinosa-Herrera}, we have the following.

\begin{prop}
For $X\in\Gamma(TM)$
and $\phi\in\Gamma(S)$, we have
\begin{eqnarray}
\sum_{i = 1}^n e_i\cdot R^{\theta,A}(X, e_i)(\phi) 
 &=& -\dfrac{1}{2} {\rm Ric}(X)
\cdot \phi + \dfrac{1}{2}\sum_{k<l}(X\lrcorner\Theta_{kl})\cdot
\kappa_{r*}(f_kf_l)\cdot\phi +{i\over 2}X\lrcorner dA\cdot \phi,\label{Ricci-curvature-identity}
\end{eqnarray}
where $\mathrm{Ric}$ denotes the Ricci tensor of $M$ and $R^{\theta,A}$ the curvature operator of the twisted spinorial connection.
\end{prop}
{\em Proof}. For $\phi=\varphi\otimes \psi$,
by \rf{curvature-twisted-spin}, 
\begin{eqnarray*}
R^{\theta,A}(X,e_\alpha)(\varphi\otimes\psi) 
   &=& 
\varphi\otimes\left[{1\over 2} \sum_{i<j} \Omega_{ij}(X,e_\alpha) e_ie_j\cdot\psi\right]\\
   &&
+\left[{1\over 2} \sum_{k<l} \Theta_{kl}(X,e_\alpha) \kappa_{r*}(f_kf_l)\cdot
\varphi\right]\otimes\psi
  +{i\over 2}\,\varphi\otimes(dA(X,e_\alpha)\psi).
\end{eqnarray*}
Multiply by $e_\alpha$ and sum over $\alpha$
\begin{eqnarray*}
\sum_\alpha e_\alpha\cdot R^{\theta,A}(X,e_\alpha)( \varphi\otimes\psi) 
   &=& 
\varphi\otimes\left[{1\over 2} \sum_\alpha \sum_{i<j} \Omega_{ij}(X,e_\alpha)  e_\alpha
e_ie_j\cdot\psi\right]\\
   &&+
{1\over 2} [\kappa_{r*}(f_kf_l)\cdot\varphi]\otimes\sum_{k<l}\left[\sum_\alpha \Theta_{kl}(X,e_\alpha)
e_\alpha \cdot \psi\right]\\
   &&
  +{i\over 2}\,\varphi\otimes \left[\sum_\alpha dA(X,e_\alpha)e_\alpha\cdot\psi\right].
\end{eqnarray*}
Now,
\begin{eqnarray*}
 {1\over 2} \sum_\alpha \sum_{i<j} \Omega_{ij}(X,e_\alpha)  e_\alpha e_ie_j  
   &=& 
  -{1\over 2} {\rm Ric}(X), \\
\kappa_{r*}(f_kf_l)\cdot\varphi \otimes{1\over 2} \sum_{k<l}\left[\sum_\alpha \Theta_{kl}(X,e_\alpha) e_\alpha \cdot \psi\right]
,
   &=& 
{1\over 2} \sum_{k<l} (X\lrcorner\Theta_{kl}) \cdot \kappa_{r*}(f_kf_l)\cdot (\varphi\otimes \psi),\\
  {i\over 2}\,\varphi\otimes \left[\sum_\alpha dA(X,e_\alpha)e_\alpha\cdot\psi\right]
   &=&
  {i\over 2}\,\varphi\otimes X\lrcorner dA\cdot\psi.
\end{eqnarray*}
\qd

\begin{prop}\label{prop: scalar curvature identity}
Let $\phi\in\Gamma(S)$. Then
\begin{equation*}
\sum_{i,j} e_ie_j\cdot R^{\theta,A}(e_i, e_j)(\phi) 
= {{\rm R}\over 2}\phi  + \sum_{k<l}\Theta_{kl}\cdot
\kappa_{r*}(f_kf_l)\cdot\phi + i\,dA\cdot\phi,
\end{equation*}
where
$\Theta_{kl}=\sum_{i<j}\Theta_{kl}(e_i,e_j)e_i\wedge e_j$ and $\mathrm{R}$ is the scalar curvature of $M$.
\end{prop}
{\em Proof}.
By \rf{Ricci-curvature-identity},
\[
\sum_{j = 1}^n e_j\cdot R^{\theta,A}(e_i, e_j)(\phi) = -\dfrac{1}{2} {\rm Ric}(e_i)
\cdot \phi + \dfrac{1}{2}\sum_j\sum_{k<l}\Theta_{kl}(e_i,e_j)e_j\cdot
\kappa_{r*}(f_kf_l)\cdot\phi 
+ {i\over 2}\, e_i\lrcorner dA\cdot\phi
.
\]
Multiplying with $e_i$ and summing over $i$, we get
\begin{eqnarray*}
\sum_{i,j} e_ie_j\cdot R^{\theta,A}(e_i, e_j)(\phi) 
&=&
 -\dfrac{1}{2} \sum_ie_i\cdot{\rm Ric}(e_i)
\cdot \phi + \dfrac{1}{2}\sum_{k<l}\left[\sum_{i,j}\Theta_{kl}(e_i,e_j)e_ie_j\right]\cdot
\kappa_{r*}(f_kf_l)\cdot\phi\\
   &&
 + {i\over 2}\, \sum_i e_i\cdot e_i\lrcorner dA\cdot\phi
.
\end{eqnarray*}
Now,
\begin{eqnarray*}
- \sum_i e_i\cdot {\rm Ric}(e_i) 
 &=& {\rm R},
\end{eqnarray*}
where ${\rm R}$ denotes the scalar curvature of $M$. 
For $k$ and $l$ fixed, 
\begin{eqnarray*}
 \sum_{i,j}\Theta_{kl}(e_i,e_j)e_ie_j 
&=& 2\sum_{i<j}\Theta_{kl}(e_i,e_j)e_ie_j\\
&=& 2\Theta_{kl},\\
  {i\over 2}\, \sum_i e_i\cdot e_i\lrcorner dA\cdot\psi
   &=&
 {i\over 2}\, \sum_{i,\alpha}  dA(e_i,e_\alpha)e_i\cdot e_\alpha\cdot\psi\\
   &=&
 {i}\, \sum_{i<\alpha}  dA(e_i,e_\alpha)e_i\cdot e_\alpha\cdot\psi\\
   &=&
 {i}\, dA\cdot\psi.
\end{eqnarray*}
\qd

Let us define
\begin{eqnarray*}
\Theta&=&\sum_{k<l}\Theta_{kl}\otimes f_kf_l \in \ext^2T^*M\otimes \ext^2 F ,\\ 
\hat{\Theta}&=&\sum_{k<l}\hat{\Theta}_{kl}\otimes f_kf_l \in \End^-(TM)\otimes \ext^2F,\\
\eta^\phi&=&\sum_{k<l}\eta_{kl}^\phi\otimes f_kf_l \in \ext^2T^*M\otimes \ext^2 F ,\\ 
\hat{\eta}^\phi&=&\sum_{k<l}\hat{\eta}_{kl}^\phi\otimes f_kf_l \in \End^-(TM)\otimes \ext^2F,
\end{eqnarray*}
where $\hat{\Theta}_{kl}$ denotes the skew-symmetric endomorphism associated to $\Theta_{kl}$ via the metric.
Denote by
\begin{eqnarray*}
\tilde\Theta
&=& (\mu_n\otimes \kappa_{r*}) (\Theta) ,
\end{eqnarray*}
the corresponding operator on twisted spinor fields. 
In order to simplify notation, we also define
\begin{eqnarray*}
 \left<{\Theta},{\eta}^\phi\right>_0 
 &=&
 \sum_{k<l}\sum_{i<j} \Theta_{kl}(e_i,e_j)\eta_{kl}^\phi(e_i,e_j),\\
 \left<\hat{\Theta},\hat{\eta}^\phi\right>_1 
 &=&
 \sum_{k<l}\tr( \hat{\Theta}_{kl}(\hat{\eta}_{kl}^\phi)^T).
\end{eqnarray*}

\begin{theo}[Twisted Schr\"odinger-Lichnerowicz Formula]\label{theo-SL-formula}
Let $\phi\in\Gamma(S)$. Then
\begin{equation}
 \dirac^{\theta,A}(\dirac^{\theta,A}(\phi)) =  \Delta^{\theta,A}(\phi)+ \dfrac{\rm R}{4}\phi +
\dfrac{1}{2}\tilde\Theta\cdot\phi +{i\over 2}dA\cdot\phi\label{SL-formula}
\end{equation}
where ${\rm R}$ is the scalar curvature of the Riemannian manifold $M$.
\end{theo}
{\em Proof}. 
Consider the
difference
\begin{eqnarray*}
\dirac^{\theta,A}(\dirac^{\theta,A}(\phi))  - \Delta^{\theta,A}(\phi)
&=& \sum_{i}\sum_{j\not =k} \left<\nabla_{e_i}e_j,e_k\right>e_ie_k\cdot\nabla_{e_i}^{\theta,A}\phi
+\sum_{i\not=j}e_ie_j\cdot\nabla_{e_j}^{\theta,A}\nabla_{e_j}^{\theta,A}\phi),
\end{eqnarray*}
since
\[\sum_j\sum_{i=k}\left<\nabla_{e_i}e_j,e_k\right>e_ie_k \nabla_{e_j}^{\theta,A}\phi =
-\sum_j\diver(e_j)\nabla_{e_j}^{\theta,A}\phi.\]
Thus,
\begin{eqnarray*}
\dirac^{\theta,A}(\dirac^{\theta,A}(\phi))  - \Delta^{\theta,A}(\phi)
&=& \sum_{j}\sum_{i<k} \left<e_j,[e_k,e_i]\right>e_ie_k\cdot\nabla_{e_i}^{\theta,A}\phi
+\sum_{i<j}e_ie_j\cdot(\nabla_{e_i}^{\theta,A}\nabla_{e_j}^{\theta,A}
-\nabla_{e_j}^{\theta,A}\nabla_{e_i}^{\theta,A})\phi\\
&=&{1\over 2}\sum_{i,j}e_ie_jR^{\theta,A}(e_i,e_j)\phi.
\end{eqnarray*}
The result follows from Proposition \ref{prop: scalar curvature identity}.
\qd

\subsection{Bochner-type arguments}

In this subsection we will prove some corollaries of the Schr\"odinger-Lichnerowicz-type formula and
Bochner type arguments (cf. \cite{Friedrich}).
For the rest of the section, let us assume that the $n$-dimensional Riemannian Spin$^{c,r}$ manifold $M$ is compact
(without border) and connected.

\subsubsection{Harmonic spinors}

A twisted spinor field $\phi\in\Gamma(S)$ such that 
\[\dirac^{\theta,A}\phi=0\]
will be called a {\em harmonic spinor}.

\begin{corol}\label{corol-no-harmonic-spinors} 
If ${\rm R}\geq 2|\tilde\Theta| + 2|dA|$ everywhere (in pointwise operator norm), then a harmonic spinor is
parallel. 
Furthermore, if the inequality is
strict at a point, then there are no non-trivial harmonic spinors
\[\ker(\dirac^{\theta,A}) =\{0\}.\]
\end{corol}
{\em Proof}. If $\phi\not=0$ is a solution of 
\[\dirac^{\theta,A}(\phi) =0,\]
by the twisted Schr\"odinger-Lichnerowicz formula \rf{SL-formula}
\[0 =   \Delta^{\theta,A}(\phi)+ \dfrac{\rm R}{4}\phi +
\dfrac{1}{2}\tilde\Theta\cdot\phi + {i\over 2}dA\cdot\phi.\]
By taking  hermitian product with $\phi$ and integrating over $M$ we get
\begin{eqnarray*}
0
&\geq&  
\int_M| \nabla^{\theta,A}\phi|^2+ \dfrac{1}{4}\int_M\left({\rm R} -
2|\tilde\Theta|-2|dA|\right)|\phi|^2.
\end{eqnarray*}
Since
\[{\rm R} -
2|\tilde\Theta|-2|dA|\geq 0,\]
then 
\[|\nabla^{\theta,A}\phi|=0,\]
so that $\phi$ is parallel,
has non-zero constant length and 
no zeroes.

Now, if
\[{\rm R} -
2|\tilde\Theta|-2|dA| > 0\]
at some point,
\[0
\geq  
 |\phi|^2\int_M\left({\rm R} -
2|\tilde\Theta|-2|dA|\right)>0.
\]
\qd

\vspace{.1in}

Now notice that
\begin{eqnarray*}
\left<\tilde\Theta\cdot\phi,\phi\right> 
&=& \left<\sum_{k<l} \left[\sum_{i<j}\Theta_{kl}(e_i,e_j)e_ie_j\right]\cdot\kappa_{r*}(f_kf_l)
\cdot\phi,\phi\right>\nonumber\\
&=& \sum_{k<l} \sum_{i<j}\Theta_{kl}(e_i,e_j)\eta_{kl}^\phi(e_i,e_j)\nonumber\\
&=& \left<{\Theta},{\eta}^{\phi}\right>_0,
\end{eqnarray*}
which is a real number dependent on the curvature of the connection $\theta$ and  the specific
spinor $\phi$.

\begin{corol}
% \label{not-massless-Dirac-spinor}
If $\phi$ is such that
\[{\rm R}|\phi|^2+ 2\left<{\Theta},{\eta}^\phi\right>_0 + 2i\left<dA\cdot\phi,\phi\right>\geq 0\] 
everywhere, and the inequality
is strict at a
point, then 
\[\dirac^\theta(\phi) \not=0.\]
\end{corol}
{\em Proof}. Suppose $\phi\not=0$ is such that 
\[\dirac^\theta(\phi) =0.\]
Then, by \rf{SL-formula}
\begin{eqnarray*}
0
&=&  
\int_M| \nabla^{\theta}\phi|^2+ \dfrac{1}{4}\int_M \left({\rm R}|\phi|^2+
2\left<{\Theta},{\eta}^\phi\right>_0+ 2i\left<dA\cdot\phi,\phi\right>\right) \geq 0,
\end{eqnarray*}
so that $\phi$ is parallel, has non-zero constant length and no zeroes.
Since
\[{\rm R}|\phi|^2+ 2\left<{\Theta},{\eta}^\phi\right>_0+ 2i\left<dA\cdot\phi,\phi\right>> 0\] 
at some point, 
\[0\geq  \int_M \left({\rm R}|\phi|^2+
2\left<{\Theta},{\eta}^\phi\right>_0+ 2i\left<dA\cdot\phi,\phi\right>\right)>0.\]
\qd

\subsubsection{Killing spinors}

A twisted spinor field $\phi\in\Gamma(S)$ is called a {\em Killing spinor} if
\[\nabla_X^{\theta,A}\phi=\mu\,X\cdot\phi\]
for all $X\in\Gamma(TM)$, and $\mu$ a complex constant.

\begin{corol} 
Suppose $\phi\not=0$ is a Killing spinor with Killing constant $\mu$.
Then $\mu$ is either real or imaginary, and
\[
\mu^2
\geq
{1\over 4n^2}\min_M({\rm R} -
2|\tilde\Theta|-2|dA|). 
\]
If the inequality is attained, then $\phi$ is parallel, i.e. $\mu=0$. 
\end{corol}
{\em Proof}. Recall that
\begin{eqnarray*}
 \dirac^{\theta,A} (\phi) 
 &=& \sum_{i=1}^n e_i\cdot\nabla_{e_i}^{\theta,A} \phi \\
&=& -n\mu\, \phi .
\end{eqnarray*}
Then, by the twisted Schr\"odinger-Lichnerowicz formula \rf{SL-formula}
\[n^2\mu^2 \phi =   \Delta^{\theta,A}(\phi)+ \dfrac{\rm R}{4}\phi +
\dfrac{1}{2}\tilde\Theta\cdot\phi + {i\over 2}dA\cdot \phi.\]
By taking  hermitian product with $\phi$ and integrating over $M$ we get
\begin{eqnarray*}
n^2\mu^2\int_M|\phi|^2
&=&  
\int_M| \nabla^{\theta,A}\phi|^2+ \int_M\dfrac{\rm  R}{4}|\phi|^2 +
\int_M\dfrac{1}{2}\left<\tilde\Theta\cdot\phi,\phi\right> + {i\over 2}\int_M\left<dA\cdot\phi,\phi\right>\\
&\geq& {1\over4}\min_M({\rm R} -
2|\tilde\Theta|-2|dA|)\int_M|\phi|^2,
\end{eqnarray*}
and the inequality follows. 
Since the right hand side of the equality above is a real number, $\mu$ must be either real or
imaginary.
Now, if the inequality is attained,
\[\int_M|\nabla^{\theta,A}\phi|^2=0 \quad\quad\mbox{and}\quad\quad \nabla^{\theta,A}\phi=0.\]
\qd

\begin{corol}
Suppose $\phi\in\Gamma(S)$ is
a Dirac eigenspinor
\[\dirac^{\theta,A} \phi = \lambda\phi.\]
Then
\begin{eqnarray*}
\lambda^2
&\geq& {n\over4(n-1)}\left(\min_M({\rm R} -
2|\tilde\Theta|-2|dA|)\right).
\end{eqnarray*}
If the lower bound is non-negative and is attained,
the spinor $\phi$ is a real Killing spinor with Killing constant
\[
\mu=\pm {1\over 2}\sqrt{{1\over n(n-1)}\min_M({\rm R} -
2|\tilde\Theta|-2|dA|)}.\]
\end{corol}
{\em Proof}. 
Let $h:M\lra\mathbb{R}$ be a fixed smooth function.
Consider the following metric connection on the twisted Spin bundle
\[\nabla_X^h\phi = \nabla_X^{\theta,A}\phi+hX\cdot\phi.\]
Let
\[
\Delta^{h}(\phi) = - \sum_{i = 1}^n
\nabla^{h}_{e_i}\nabla^{h}_{e_i}\phi - \sum_{i = 1}
\diver(e_i)\nabla^{h}_{e_i}\phi, 
\]
be the Laplacian for this connection and recall that
\[|\nabla^h\phi|^2= \sum_{i=1}^n |\nabla_{e_i}^{\theta,A}\phi + he_i\cdot \phi|^2.\]
Then, by \rf{SL-formula}
\begin{eqnarray*}
 (\dirac^{\theta,A} - h)\circ (\dirac^{\theta,A} - h) (\phi)
&=&\dirac^{\theta,A}(\dirac^{\theta,A}\phi)-2h\dirac^{\theta,A}\phi-{\rm grad}(h)\cdot\phi+h^2\phi\\
&=&\Delta^{\theta,A}(\phi)+ \dfrac{\rm R}{4}\phi +
\dfrac{1}{2}\tilde\Theta\cdot\phi+{i\over 2}dA\cdot\phi
-2h\dirac^{\theta,A}\phi-{\rm grad}(h)\cdot\phi+h^2\phi.
\end{eqnarray*}
On the other hand,
\[\Delta^h\phi =  \Delta^{\theta,A}\phi -2h\dirac^{\theta,A}\phi-{\rm grad}(h)\cdot\phi+nh^2\phi.\]
Thus
\[ 
(\dirac^{\theta,A} - h)\circ (\dirac^{\theta,A} - h) (\phi)
=\Delta^{h}(\phi)+ \dfrac{\rm R}{4}\phi +
\dfrac{1}{2}\tilde\Theta\cdot\phi+ {i\over 2}dA\cdot\phi
+(1-n)h^2\phi 
\]
By using $\dirac^{\theta,A}\phi =\lambda\phi$, setting $h={\lambda\over n}$, taking hermitian product
with $\phi$ and integrating over $M$ we get
\[\lambda^2\left({n-1\over n}\right)^2\int_M|\phi|^2=\int_M|\nabla^{\lambda/n}\phi|^2+\lambda^2
{1-n\over n^2}\int_M|\phi|^2+\int_M\dfrac{\rm R}{4}|\phi|^2 +
\int_M\dfrac{1}{2}\left<\tilde\Theta\cdot\phi,\phi\right>
+{i\over 2}\int_M \left<dA\cdot\phi,\phi\right>\]
so that
\begin{eqnarray*}
\lambda^2\left({n-1\over n}\right)\int_M|\phi|^2
&\geq& {1\over4}\min_M({\rm R} -
2|\tilde\Theta|-2|dA|)\int_M|\phi|^2.
\end{eqnarray*}

If the lower bound is attained,
\[\int_M|\nabla^{\lambda/n}\phi|^2= 0,
\]
i.e.
\[\nabla^{\lambda/n}\phi= 0.
\]
\qd

\section{CR structures of arbitrary codimension}\label{sec: CR structures}

In this section we will explore the twisted spinorial geometry associated to almost CR structures. 
We carry out the spinorial characterization and explore some
integrability conditions of almost CR structures implied by assuming the typical conditions on spinors,
such as being parallel or Killing, but just in
prescribed directions. 

\subsection{Spinorial characterization of almost CR (hermitian) structures}\label{subsec: spinorial
characterization almost CR hermitian}

\begin{defi} Let $M$ be a smooth $(2m+r)$-dimensional smooth manifold.
\begin{itemize}
\item 
An {\em almost CR structure} on a manifold $M$ consists of
a sub-bundle $D\subset TM$ and
a bundle automorphism $J$ of $D$ such that $J^2=-{\rm Id}_D$.

\item An {\em almost CR hermitian structure} on $M$ is an almost CR structure whose almost complex structure is
orthogonal with respect to the metric.
\end{itemize}
\end{defi}

{\bf Remark}. Given an almost CR structure on $M$ we can introduce an (auxiliary) compatible
metric as follows. Take any Riemannian metric $g_0$ on $M$ and consider the
orthogonal complement $D^{\perp}$ of $D$ with respect to this metric.
Let $g_1$ and $g_2$ denote the restrictions of $g_0$ to $D$ and $D^{\perp}$
respectively. Average $g_1$ with respect to $J$ and call it $g_3$. Finally,
consider the metric $g=g_3\oplus g_2$.

\begin{defi} 
Let $M$ be an oriented Riemannian Spin$^{c,r}(n)$ manifold and $S$ the associated twisted spinor bundle.
A (nowhere zero) spinor field
$\phi\in \Gamma(S)$ is called {\em partially pure} if $\phi_x\in S_x$ is partially
pure at each point $x\in M$.
\end{defi}

\begin{theo}\label{theo: characterization almost CR}
 Let $M$ be an oriented $n$-dimensional Riemannian manifold. Then the following two statements are
equivalent:
\begin{itemize}
 \item[(a)] $M$ admits a twisted Spin$^{c,r}$ structure carrying a partially pure spinor field
$\phi\in
\Gamma(S)$, where $S$ denotes the associated tiwsted spinor bundle.
 \item[(b)] $M$ admits an almost CR hermitian structure of codimension $r$.
\end{itemize}
\end{theo}

{\em Proof}. 
If the manifold $M$ admits a partially pure spinor field $\phi\in
\Gamma(S)$, the subspaces
$V_x^\phi$ 
determine a smooth distribution of even rank $n-r$ carrying an almost complex structure. 

Conversely, if $M$ has an orthogonal almost CR hermitian structure of codimension $r$, the tangent bundle
decomposes
orthogonally as
\[TM=D\oplus D^\perp,\]
where $D$ has real rank $2m=n-r$ and admits an almost complex structure, and
$D^\perp$ is the oriented orthogonal complement.
The structure group of the Riemannian manifold $M$ reduces to $U(m)\times SO(r)$ and,
by Lemma \ref{lemma:subgroup2}, there is a monomorphism
\[U(m)\times SO(r)\hookrightarrow Spin^{c,r}(2m+r)\]
with image $\widehat{U(m)\times SO(r)}$,
which allows us to associate a $Spin^{c,r}(n)$ principal bundle $P$ on $M$, i.e. a Spin$^{c,r}$
structure.
Note that the corresponding twisted spinor bundle $S$ decomposes under $\widehat{U(m)\times SO(r)}$ as follows
\begin{eqnarray*}
 S&=&
\left[\kappa_D^{-1/2}\otimes \Sigma(D^\perp)\right]\otimes \Delta(M) \\
   &=&
\left[\kappa_D^{-1/2}\otimes \Sigma(D^\perp)\right]\otimes\Delta(D^\perp)\otimes \Delta(D)\\
   &=&
\left[\kappa_D^{-1/2}\otimes \Sigma(D^\perp)\right]\otimes\Delta(D^\perp)\otimes \left[(\ext^*
D^{0,1})\otimes \kappa_D^{1/2}\right]\\
   &=&
\left[ \Sigma(D^\perp)\otimes\Delta(D^\perp)\right]\otimes \left[\ext^*D^{0,1}\right],
\end{eqnarray*}
where $\kappa_D=\ext^mD^{1,0}$.
We see that it contains a rank 1 trivial subbundle
generated by the partially pure spinor given in  
\rf{eq: canonical partilly pure spinor}
with stabilizer $\widehat{U(m)\times SO(r)}$, i.e. $M$ admits a global
partially pure spinor field.
\qd

{\bf Example}.
Recall from Subsection \ref{subsec: certain homogeneous spaces} that
\[T_{{\rm Id}}\mathcal{G}_{m,1,r} \cong  
 [[\ext^2\mathbb{C}^m]] \oplus [[\mathbb{C}^m]]\otimes \mathbb{R}^r 
\oplus [[\mathbb{C}^m]]
\oplus
\mathbb{R}^r
\]
For the sake of clarity, consider $m=2$, $r=2$ and
$\mathbb{R}^7=\mathbb{R}^4\oplus\mathbb{R}^2\oplus\mathbb{R}^1$,
where the first summand $\mathbb{R}^4$ is endowed with the standard complex structure
\[\left(
\begin{array}{cccc}
0 & -1 &  & \\
1 & 0 &  & \\
 &  & 0 & -1\\
 &  & 1 & 0
        \end{array}
\right).\]
The different summands in the decomposition
\[T_{{\rm Id}}\mathcal{G} \cong  
 [[\ext^2\mathbb{C}^2]] \oplus [[\mathbb{C}^2]]\otimes \mathbb{R}^2 
\oplus [[\mathbb{C}^2]]
\oplus
\mathbb{R}^2.
\]
correspond to skew-symmetric matrices as follows:
\begin{eqnarray*}
[[\ext^2\mathbb{C}^2]] 
   &=&
  \left\{\left(\begin{array}{ccccccc}
0 & 0 & b_1 & b_2 &  &  & \\
 & 0 & b_2 & -b_1 &  &  & \\
-b_1 & -b_2 & 0 &  &  &  & \\
-b_2 & b_1 &  & 0 &  &  & \\
 &  &  &  & 0 &  & \\
 &  &  &  &  & 0 & \\
 &  &  &  &  &  & 0
                                       \end{array}
\right): \,\, b_1,b_2\in\mathbb{R}\right\},  \\ 
{}[[\mathbb{C}^2]]\otimes \mathbb{R}^2
   &=&
\left\{\left(\begin{array}{ccccccc}
0 &  &  &  & c_1 & c_2 & \\
 & 0 &  &  & c_3 & c_4 & \\
 &  & 0 &  & c_5 & c_6 & \\
 &  &  & 0 & c_7 & c_8 & \\
-c_1 & -c_3 & -c_5 & -c_7 & 0 &  & \\
-c_2 & -c_4 & -c_6 & -c_8 &  & 0 & \\
 &  &  &  &  &  & 0
             \end{array}
\right):\,\, c_j\in\mathbb{R},\,\, j=1,\ldots,8\right\},\\
{}[[\mathbb{C}^2]]
   &=&
 \left\{\left(\begin{array}{ccccccc}
0 &  &  &  &  &  & d_1\\
 & 0 &  &  &  &  & d_2\\
 &  & 0 &  &  &  & d_3\\
 &  &  & 0 &  &  & d_4\\
 &  &  &  & 0 &  & \\
 &  &  &  &  & 0 & \\
-d_1 & -d_2 & -d_3 & -d_4 &  &  & 0
              \end{array}
\right): \,\, d_j\in\mathbb{R},\, j=1,\ldots,4\right\},\\
\mathbb{R}^2
   &=&
 \left\{\left(\begin{array}{ccccccc}
0 &  &  &  &  &  & \\
 & 0 &  &  &  &  & \\
 &  & 0 &  &  &  & \\
 &  &  & 0 &  &  & \\
 &  &  &  & 0 &  & \delta_1\\
 &  &  &  &  & 0 & \delta_2\\
 &  &  &  & -\delta_1 & -\delta_2 & 0
              \end{array}
\right):\,\, \delta_1,\delta_2\in\mathbb{R}\right\}.
\end{eqnarray*}
The induced complex structure on $[[\ext^2\mathbb{C}^2]] \oplus [[\mathbb{C}^2]]\otimes \mathbb{R}^2 
\oplus [[\mathbb{C}^2]]$, which respects each summand, is
\[J\left(\begin{array}{ccccccc}
0 & 0 & b_1 & b_2 & c_1 & c_2 & d_1\\
 & 0 & b_2 & -b_1 & c_3 & c_4 & d_2\\
 &  & 0 & 0 & c_5 & c_6 & d_3\\
 &  &  & 0 & c_7 & c_8 & d_4\\
 &  &  &  & 0 & 0 & 0 \\
 &  &  &  &  & 0 & 0\\
 &  &  &  &  &  & 0
         \end{array}
\right)
=
\left(
\begin{array}{ccccccc}
0 & 0 & -b_2 & b_1 & -c_3 & -c_4 & -d_2\\
 & 0 & b_1 & b_2 & c_1 & c_2 & d_1\\
 &  & 0 & 0 & -c_7 & -c_8 & -d_4\\
 &  &  & 0 & c_5 & c_6 & d_3\\
 &  &  &  & 0 & 0 & 0\\
 &  &  &  &  & 0 & 0\\
 &  &  &  &  &  & 0
      \end{array}
\right)
\]
where we have only written the upper triangle part for notational simplicity. 

Thus, this example gives us several candidates of distributions carrying an almost complex structure, 
as well as their orthogonal complements (with respect to the natural metric):
\begin{eqnarray*}
&&\left\{\begin{array}{lll}
D_1 & = & [[\ext^2\mathbb{C}^m]] \\
D_1^\perp & = & [[\mathbb{C}^m]]\otimes \mathbb{R}^r \oplus [[\mathbb{C}^m]] \oplus \mathbb{R}^r
         \end{array}
 \right.\\
&&\left\{\begin{array}{lll}
D_2 & = & [[\mathbb{C}^m]]\otimes \mathbb{R}^r \\
D_2^\perp & = & [[\ext^2\mathbb{C}^m]] \oplus [[\mathbb{C}^m]] \oplus \mathbb{R}^r
         \end{array}
 \right.\\
&&\left\{\begin{array}{lll}
D_3 & = & [[\mathbb{C}^m]] \\
D_3^\perp & = & [[\ext^2\mathbb{C}^m]] \oplus [[\mathbb{C}^m]]\otimes \mathbb{R}^r \oplus \mathbb{R}^r
         \end{array}
 \right.\\
&&\left\{\begin{array}{lll}
D_4 & = & [[\ext^2\mathbb{C}^m]] \oplus [[\mathbb{C}^m]]\otimes \mathbb{R}^r \\
D_4^\perp & = & [[\mathbb{C}^m]] \oplus \mathbb{R}^r
         \end{array}
 \right.\\
&&\left\{\begin{array}{lll}
D_5 & = & [[\ext^2\mathbb{C}^m]] \oplus [[\mathbb{C}^m]] \\
D_5^\perp & = & [[\mathbb{C}^m]]\otimes \mathbb{R}^r \oplus \mathbb{R}^r
         \end{array}
 \right.\\
&&\left\{\begin{array}{lll}
D_6 & = & [[\mathbb{C}^m]]\otimes \mathbb{R}^r \oplus [[\mathbb{C}^m]] \\
D_6^\perp & = & [[\ext^2\mathbb{C}^m]] \oplus \mathbb{R}^r
         \end{array}
 \right.\\
&&\left\{\begin{array}{lll}
D_7 & = & [[\ext^2\mathbb{C}^m]] \oplus [[\mathbb{C}^m]]\otimes \mathbb{R}^r \oplus [[\mathbb{C}^m]] \\
D_7^\perp & = &  \mathbb{R}^r
         \end{array}
 \right.
\end{eqnarray*}

By computing the Lie brackets at the Lie algebra level, we see that the distributions $D_1, D_4, D_6^\perp$ 
and $D_7^\perp$ are involutive with their foliations corresponding to the fibers of the following four fibrations
\[
\begin{array}{ccc}
{SO(2m)\over U(m)} & \hookrightarrow & {SO(2m+r+1)\over U(m)\times SO(r)}\\
 &  & \downarrow\\
 &  & {SO(2m+r+1)\over SO(2m)\times SO(r)},\\
&&\\
{SO(2m+r)\over U(m)\times SO(r)} & \hookrightarrow & {SO(2m+r+1)\over U(m)\times SO(r)}\\
 &  & \downarrow\\
 &  & S^{2m+r},\\
&&\\
{SO(2m)\over U(m)}\times S^r & \hookrightarrow & {SO(2m+r+1)\over U(m)\times SO(r)}\\
 &  & \downarrow\\
 &  & {SO(2m+r+1)\over SO(2m)\times SO(r+1)},\\
&&\\
S^r & \hookrightarrow & {SO(2m+r+1)\over U(m)\times SO(r)}\\
 &  & \downarrow\\
 &  & {SO(2m+r+1)\over U(m)\times SO(r+1)},
\end{array}
\]
respectively.

\subsection{Adapted connection for almost CR-hermitian manifolds}\label{subsec: adapted connections}

Before we proceed with the characterizations of integrability conditions, we need to give (at least) a choice of connection on the relevant bundles 
of an almost CR hermitian manifold, or equivalently, on a Spin$^{c,r}$ manifold carrying a partially pure spinor.

As we mentioned earlier,
we can adapt a
metric on an almost CR manifold $M$ in order to make it an almost CR hermitian manifold.
Let us fix one such metric and its Levi-Civita connection 1-form $\omega$ and covariant derivative $\nabla$.
The metric determines the orthogonal complement $D^\perp$ and gives us a covariant derivative as follows
\begin{eqnarray*}
\nabla_X^{D^\perp} : \Gamma(D^\perp)&\longrightarrow& \Gamma( D^\perp)\\
  W &\mapsto& {\rm proj}_{D^\perp}(\nabla_XW)  
\end{eqnarray*}
for $X\in\Gamma(TM)$, whose local connection 1-forms and curvature 2-forms will be denoted by 
$\theta_{kl}^{D^\perp}$ and $\Theta_{kl}^{D^\perp}$ respectively, $1\leq k<l\leq r$.
The analogous connection on $D$ is given by the covariant derivative
\begin{eqnarray*}
\nabla_X^{D} : \Gamma(D)&\longrightarrow& \Gamma( D)\\
  W &\mapsto& {\rm proj}_{D}(\nabla_XW)  
\end{eqnarray*}
for $X\in\Gamma(TM)$.
However, we need to induce a connection on $\kappa_D^{-1}$.
Thus, we consider the hermitian connection for $(D,\left<,\right>,J)$ defined by
\begin{eqnarray*}
\tilde\nabla^D_XY 
   &=& 
  \nabla^D_XY +{1\over 2} (\nabla^D_XJ)(JY). 
\end{eqnarray*}
so that
\[\tilde\nabla^D J=0.\]
$\tilde\nabla^D$ induces a covariant derivative $\tilde\nabla^{\kappa_D^{-1}}$ on the anticanonical bundle $\kappa_D^{-1}$ of
$D$,
whose local connection 1-form will be denoted by $i\tilde{A}^D$.
More precisely, if $(e_1,\ldots,e_n)$ is a local orthonormal frame of $TM$ such
that
\begin{eqnarray*}
D
   &=&
  \span(e_1,\ldots, e_{2m}),\\
e_{2s}
   &=& J (e_{2s-1}),\\
D^\perp
   &=&
  \span(e_{2m+1},\ldots, e_{2m+r}),
\end{eqnarray*}
for $1\leq s\leq m$ and $1\leq k<l\leq r$. 
and the matrix of connection 1-forms of $\tilde\nabla^D$ is 
\[
\left(
\begin{array}{ccccc}
0 & \tilde{\omega}_{1,2} &  & \tilde{\omega}_{1,2m-1} & \tilde{\omega}_{1,2m} \\
-\tilde{\omega}_{12} & 0 &  & -\tilde{\omega}_{1,2m} & \tilde{\omega}_{1,2m-1} \\
 &  & \ddots &  &  \\
-\tilde{\omega}_{1,2m-1} & \tilde{\omega}_{1,2m} &  & 0 & \tilde{\omega}_{2m-1,2m} \\
-\tilde{\omega}_{1,2m} & -\tilde{\omega}_{1,2m-1} & & -\tilde{\omega}_{2m-1,2m} & 0
\end{array}
\right),
\] 
the induced connection on $\kappa_D^{-1}=\ext^m D^{0,1}$ is
\[i\tilde A=-i[\tilde{\omega}_{1,2}+\cdots+\tilde{\omega}_{2m-1,2m}].
\]
By using $\nabla$, $\nabla^{D^\perp}$ and the unitary connection $i\tilde{A}^D$, we can define a
connection $\nabla^{S}$
on the 
globally defined twisted spinor vector bundle
$S=\left[\kappa_D^{-1/2}\otimes \Sigma(D^\perp)\right]\otimes \Delta(M) $ which is compatible with
Clifford
multiplication.

\subsection{Spinorial characterization of integrability}

\begin{defi} Let $M$ be a smooth $2m+r$ dimensional smooth manifold.
An almost CR structure is called a {\em CR structure} if
for every $X,Y\in \Gamma (D)$
\begin{itemize}
 \item $[X,Y]-[J(X),J(Y)]\in \Gamma(D)$,
 \item $[J(X),Y]+[X,J(Y)]\in \Gamma(D)$,
 \item $J([X,Y]-[J(X),J(Y)])=[J(X),Y]+[X,J(Y)]$.
\end{itemize}
\end{defi}

{\bf Example}. By computing the relevant combinations of brackets one can check that the 
distributions $D_1, D_4, D_5$ and $D_7$  on $\mathcal{G}_{m,1,r}$ are CR-integrable.

\begin{theo}\label{theo: integrable} 
Let $M$ be an oriented $n$-dimensional Riemannian manifold. The following are equivalent:
\begin{enumerate}
 \item[{\rm (i)}] $M$ is endowed with a
CR hermitian structure of codimension $r$.
\item[{\rm (ii)}] $M$ admits a twisted Spin$^{c,r}(n)$ structure and a twisted spinor bundle $S$ 
carrying a partially pure spinor field $\phi\in\Gamma(S)$ which satisfies  
\[(X-iJ^\phi(X))\cdot\nabla_{(Y-iJ^\phi(Y))}^S\phi = (Y-iJ^\phi(Y))\cdot\nabla_{(X-iJ^\phi(X))}^S\phi,\]
for every $X, Y \in \Gamma(V^\phi)$, where $\nabla^S$ is the covariant drivative described in subsection \ref{subsec: adapted connections}.
\end{enumerate}
\end{theo}
{\em Proof.} 
First, let us assume (i), i.e. $M$ admits a CR hermitian structure. 
By Theorem \ref{theo: characterization almost CR}, 
$M$ admits a twisted spinor vector bundle $S=\left[\kappa_D^{-1/2}\otimes \Sigma(D^\perp)\right]\otimes \Delta(M) $  
carrying a
partially pure spinor
field $\phi\in\Gamma(S)$ such that $V^\phi=D$, $J^\phi=J$ 
and
\[(X-iJX)\cdot\phi =0\] 
for every $X\in\Gamma(D)$. 
By differentiating this identity 
\begin{equation}
\left(\nabla_Y X  -i\nabla_Y(JX)\right)\cdot\phi + \left(X - i JX\right)\cdot\nabla_Y^{S}\phi 
=0,\label{eq:1}
\end{equation}
and similarly
\begin{equation}
\left(\nabla_X Y  -i\nabla_X(JY)\right)\cdot\phi + \left(Y - i JY\right)\cdot\nabla_X^{S}\phi 
=0,\label{eq:2}
\end{equation}
By subtracting \rf{eq:1} from \rf{eq:2} 
\begin{equation}
\left([X,Y]  -i\nabla_X(JY)+i\nabla_Y(JX)\right)\cdot\phi = 
 \left(X - i JX\right)\cdot\nabla_Y^{S}\phi -\left(Y - i JY\right)\cdot\nabla_X^{S}\phi 
, \label{eq:3} 
\end{equation}
By substituting $X$ with $JX$, and $Y$ with $JY$ in \rf{eq:3}
\begin{equation}
\left([JX,JY]  +i\nabla_{JX}(Y)-i\nabla_{JY}(X)\right)\cdot\phi = 
 -\left(X - iJX\right)\cdot\nabla_{-iJY}^{S}\phi +\left(Y - i JY\right)\cdot\nabla_{-iJX}^{S}\phi 
, \label{eq:4} 
\end{equation}
Subtract \rf{eq:4} from \rf{eq:3}
\begin{equation}
\left([X,Y] -[JX,JY] -i([X,JY]+[JX,Y])
\right)\cdot\phi=
 \left(X - iJX\right)\cdot\nabla_{Y-iJY}^{S}\phi -\left(Y - i JY\right)\cdot\nabla_{X-iJX}^{S}\phi .
\label{eq: characterization CR integrability 1}
\end{equation}
Since $[X,Y] -[JX,JY]\in\Gamma(D)$
\begin{eqnarray*}
([X,Y] -[JX,JY])\cdot \phi 
   &=& i\, J([X,Y] -[JX,JY])\cdot \phi\\
   &=& i\, ([J(X),Y]+[X,J(Y)])\cdot \phi,
\end{eqnarray*}
so that the left hand side of \rf{eq: characterization CR integrability 1} vanishes.

Conversely, let us assume (ii). Then, the subbundle $V^\phi$ together with
its endomorphism $J^\phi$ provide an almost CR hermitian structure on $M$. By considering the equation
\[(X-iJ^\phi X)\cdot\phi =0\] 
for all $X\in V^\phi$, and performing the same calculations as before, we arrive at 
\begin{eqnarray*}
\left([X,Y] -[J^{\phi} X,J^{\phi}Y] -i([X,J^{\phi}Y]+[J^{\phi}X,Y])
\right)\cdot\phi
   &=&
 \left(X - iJ^{\phi}X\right)\cdot\nabla_{Y-iJ^{\phi}Y}^{S}\phi -\left(Y - i
J^{\phi}Y\right)\cdot\nabla_{X-iJ^{\phi}X}^{S}\phi\\ 
&=&0 , 
\end{eqnarray*}
i.e.
\[([X,Y] -[J^{\phi} X,J^{\phi}Y])\cdot\phi =i([X,J^{\phi}Y]+[J^{\phi}X,Y])
\cdot\phi,\]
which implies 
\begin{itemize}
 \item $[X,Y]-[J^{\phi}(X),J^{\phi}(Y)]\in \Gamma(V^\phi)$,
 \item $[J^{\phi}(X),Y]+[X,J^{\phi}(Y)]\in \Gamma(V^\phi)$,
 \item $J^{\phi}([X,Y]-[J^{\phi}(X),J^{\phi}(Y)])=[J^{\phi}(X),Y]+[X,J^{\phi}(Y)]$,
\end{itemize}
since $\phi$ is a partially pure spinor.
\qd

\subsection{$D$-parallel partially pure spinor}\label{subsec: D parallelness}

The following theorem is motivated by the condition
\[\nabla_X^S\phi =0\]
for all $X\in\Gamma(D)$, i.e. $\phi$ being $D$-parallel.

\begin{theo}\label{theo: D-parallel}
Let $M$ be an oriented $n$-dimensional Riemannian manifold. The following are equivalent:
\begin{enumerate}
\item[{\rm (i)}] $M$ admits a twisted Spin$^{c,r}(n)$ structure and a twisted spinor bundle $S$
carrying a
partially pure spinor field $\phi\in\Gamma(S)$ satisfying
\[ (Y-iJ^\phi(Y))\cdot\nabla_{X}^S\phi=0\]
for every $X, Y \in \Gamma(V^\phi)$, where $\nabla^S$ is the covariant derivative described in subsection \ref{subsec: adapted connections}.
 \item[{\rm (ii)}] $M$ is endowed with an
almost CR hermitian structure of codimension $r$, where 
$D$ and $J$ are $D$-parallel. (In particular, $J$ restricts to a 
K\"ahler structure on each leaf of the integral foliation of $D$, and $D^\perp$ is $D$-parallel.) 
\end{enumerate}
\end{theo}
{\em Proof.}
Let us assume (i) and $D=V^\phi$, $D^\perp=(V^\phi)^\perp$, $J=J^\phi$.
Since
\[(Y-iJ^{\phi}Y)\cdot\phi =0\] 
for every $Y\in\Gamma(V^\phi)$, if 
 $X\in\Gamma(V^\phi)$
\begin{eqnarray*}
 0
   &=& \nabla_X^S((Y-iJ^{\phi}Y)\cdot \phi)\\
   &=& (\nabla_XY-i\nabla_X(J^{\phi}Y))\cdot \phi+ (Y-iJ^{\phi}Y)\cdot \nabla_X^{S}\phi\\
   &=& (\nabla_XY-i\nabla_X(J^{\phi}Y))\cdot \phi,
\end{eqnarray*}
which means
\begin{eqnarray*}
 \nabla_XY &\in& D,\\
 \nabla_X(JY)&=& J(\nabla_XY).
\end{eqnarray*}
i.e. $D$ and $J$ are $D$-parallel so that the leaves of this totally geodesic foliation
are K\"ahler manifolds. 
If $u\in\Gamma(D^\perp)$
\[\left<Y,u\right>=0\]
for every $Y\in\Gamma(D)$, so that for every $X\in\Gamma(D)$
\begin{eqnarray*}
 0
   &=& 
  X\left<Y,u\right>\\
   &=& 
  \left<Y,\nabla_Xu\right>
\end{eqnarray*}
since $D$ is $D$-parallel, thus showing that
$\nabla_Xu\in\Gamma(D^\perp)$.

Conversely, if $M$ admits an almost CR hermitian structure.
By Theorem \ref{theo: characterization almost CR}, 
$M$ admits a twisted Spin structure and a twisted spinor bundle $S$ endowed with a connection $\nabla^{S}$, carrying a partially pure spinor
field $\phi\in\Gamma(S)$ such that $V^\phi=D$, $J^\phi=J$ 
and
\[(Y-iJY)\cdot\phi =0\] 
for every $Y\in\Gamma(D)$.
Thus, for $X\in\Gamma(D)$,
\begin{eqnarray*}
 0
   &=& \nabla_X((Y-iJY)\cdot \phi)\\
   &=& (\nabla_XY-i J(\nabla_XY))\cdot \phi+ (Y-iJ(Y)\cdot \nabla_X^{S}\phi\\
   &=&  (Y-iJ(Y))\cdot \nabla_X^{S}\phi
\end{eqnarray*}
since $J$ is $D$-parallel.
As before, $D^\perp$ is $D$-parallel.
\qd

{\bf Example}. The space $\mathcal{G}_{m,1,r}$ admits the CR distribution $D_1$ satisfying the hypotheses of Theorem \ref{theo: D-parallel}, 
as can be seen from the fibration:
\[
\begin{array}{ccc}
{SO(2m)\over U(m)} & \hookrightarrow & {SO(2m+r+1)\over U(m)\times SO(r)}\\
 &  & \downarrow\\
 &  & {SO(2m+r+1)\over SO(2m)\times SO(r)}.
\end{array}
\]

\vspace{.in}

When the partially pure spinor is parallel, we can actually say more about the foliation leaves' Ricci curvature.

\begin{theo}\label{theo: Ricci tensor D-parallel spinor}
Let $M$ be a Spin$^{c,r}$ $n$-dimensional Riemannian manifold such that
its twisted spinor bundle $S$
admits a partially pure spinor field $\phi\in\Gamma(S)$ satisfying
\[ \nabla_{X}^S\phi=0\]
for every $X \in \Gamma(V^\phi)$, where $\nabla^S$ is the covariant derivative described in subsection \ref{subsec: adapted connections}.
Then 
\begin{enumerate}
 \item The Ricci tensor of $V^\phi$ satisfies
\begin{equation}
{\rm Ric}^{V^\phi} =\left[{\rm proj}_{V^\phi}\circ\widehat{dA}|_{V^\phi}\right]\circ J^\phi,\label{eq: Ricci form} 
\end{equation}
where $\widehat{dA}$ denotes the skew-symmetric endomorphism  determined by $dA$ (the curvature of the connection $1$-form on the auxiliary principal $U(1)$ bundle) and metric
dualization.

 \item The scalar curvature is given by
\[{\rm R}^{V^\phi}
=  \tr\left(\left[{\rm proj}_{V^\phi}\circ\widehat{dA}|_{V^\phi}\right]\circ J^\phi\right).
\] 
 \item If the connection $A$ on the auxiliary bundle $L$ is flat along an integral leaf of $V^\phi$, then the
leaf is Calabi-Yau. 
\end{enumerate}
\end{theo}

{\bf Remark}. The identity \rf{eq: Ricci form} tells us that ${\rm proj}_{V^\phi}\circ\widehat{dA}|_{V^\phi}$, restricted to the leaves of the corresponding foliation,
equals their Ricci form.

{\em Proof}. 
 Since $\phi$ is partially pure, $n=2m+r$ where ${\rm rank}(V^\phi)=2m$ and
${\rm rank}((V^\phi)^\perp)=r$.
Let $(e_1,\ldots, e_n)$ and $(f_1,\ldots f_r)$ be local orthonormal frames of $TM$ and $F$ respectively, such
that
\begin{eqnarray*}
V^\phi
   &=&
  \span(e_1,\ldots, e_{2m}),\\
e_{2j}
   &=& J^\phi (e_{2j-1}),\\
(V^\phi)^\perp
   &=&
  \span(e_{2m+1},\ldots, e_{2m+r}),
  \\ 
\eta_{kl}^\phi 
   &=&
  e_{2m+k}\wedge e_{2m+l},
\end{eqnarray*}
for $1\leq j\leq m$ and $1\leq k<l\leq r$.
If $X\in\Gamma(V^\phi)$ and $1\leq\alpha\leq 2m$ then, by Theorem \ref{theo: D-parallel}, $[X,e_\alpha]\in\Gamma(V^\phi)$ 
and
\begin{eqnarray*}
R^M(X,e_\alpha)e_i 
   \in 
  \Gamma(V^\phi) \quad\quad&&\mbox{if $1\leq i \leq 2m$,} \\ 
R^M(X,e_\alpha)e_i 
   \in 
  \Gamma((V^\phi)^\perp) \quad\quad&&\mbox{if $2m+1\leq i \leq 2m+r$.} 
\end{eqnarray*}
so that
\begin{eqnarray*}
\left<R^M(X,e_\alpha)e_i, e_j\right>=0 
    \quad\quad&&\mbox{if $1\leq i \leq 2m$, $2m+1\leq j \leq 2m+r$,} \\ 
\left<R^M(X,e_\alpha)e_i,e_j\right> =0
    \quad\quad&&\mbox{if $2m+1\leq i \leq 2m+r$, $1\leq j \leq 2m$.} 
\end{eqnarray*}
For $\phi$, 
\begin{eqnarray*}
0
   &=&
R^{\theta,A}(X,e_\alpha)\phi \\
   &=& 
{1\over 2} \sum_{1\leq i<j\leq n} \left<R^M(X,e_\alpha)e_i,e_j\right> e_ie_j\cdot\phi
+{1\over 2} \sum_{1\leq k<l\leq r} \Theta_{kl}(X,e_\alpha) \kappa_{r*}(f_kf_l)\cdot
\phi
  +{i\over 2}\,dA(X,e_\alpha)\phi\\
   &=& 
{1\over 2} \sum_{1\leq i<j\leq 2m} \left<R^M(X,e_\alpha)e_i,e_j\right> e_ie_j\cdot\phi
+{1\over 2} \sum_{1\leq k<l\leq r} \left<R^M(X,e_\alpha)e_{2m+k},e_{2m+l}\right> e_{2m+k}e_{2m+l}\cdot\phi\\
  &&
+{1\over 2} \sum_{1\leq k<l\leq r} \Theta_{kl}(X,e_\alpha) \kappa_{r*}(f_kf_l)\cdot
\phi
  +{i\over 2}\,dA(X,e_\alpha)\phi\\
   &=& 
{1\over 2} \sum_{1\leq i<j\leq 2m} \left<R^M(X,e_\alpha)e_i,e_j\right> e_ie_j\cdot\phi
+{1\over 2} \sum_{1\leq k<l\leq r} \left<R^M(X,e_\alpha)e_{2m+k},e_{2m+l}\right> \kappa_{r*}(f_kf_l)\cdot\phi\\
  &&
+{1\over 2} \sum_{1\leq k<l\leq r} \Theta_{kl}(X,e_\alpha) \kappa_{r*}(f_kf_l)\cdot
\phi
  +{i\over 2}\,dA(X,e_\alpha)\phi,
\end{eqnarray*}
where $\Theta_{kl}$ denote the local curvature 2-forms of the auxiliary connection on $P_{SO(r)}$.
Multiply by $e_\alpha$ and sum over $\alpha$, $1\leq\alpha\leq 2m$,
\begin{eqnarray*}
0
   &=&
 \sum_{\alpha=1}^{2m}\sum_{1\leq i<j\leq 2m} \left<R^M(X,e_\alpha)e_i,e_j\right> e_\alpha e_ie_j\cdot\phi
+ \sum_{\alpha=1}^{2m}\sum_{1\leq k<l\leq r} \left<R^M(X,e_\alpha) e_{2m+k},e_{2m+l}\right>e_\alpha\cdot
\kappa_{r*}(f_kf_l)\cdot\phi\\
  &&
+ \sum_{\alpha=1}^{2m}\sum_{1\leq k<l\leq r} \Theta_{kl}(X,e_\alpha) e_\alpha\cdot\kappa_{r*}(f_kf_l)\cdot
\phi
  +{i}\sum_{\alpha=1}^{2m}\,dA(X,e_\alpha)e_\alpha\cdot\phi\\
   &=&
-{\rm Ric}^{V^\phi}(X)\cdot\phi
+ \sum_{\alpha=1}^{2m}\sum_{1\leq k<l\leq r} \left<R^M(X,e_\alpha) e_{2m+k},e_{2m+l}\right>e_\alpha\cdot
\kappa_{r*}(f_kf_l)\cdot\phi\\
  &&
+ \sum_{\alpha=1}^{2m}\sum_{1\leq k<l\leq r} \Theta_{kl}(X,e_\alpha) e_\alpha\cdot\kappa_{r*}(f_kf_l)\cdot
\phi
  +{i}\sum_{\alpha=1}^{2m}\,dA(X,e_\alpha)e_\alpha\cdot\phi.
\end{eqnarray*}
By taking the real part of the hermitian inner product with $e_i\cdot\phi$, $1\leq i\leq 2m$,
\begin{eqnarray*}
{\rm Re}\left<{\rm Ric}^{V^\phi}(e_j)\cdot\phi,e_i\cdot\phi\right> 
&=& \left<{\rm Ric}^{V^\phi}(e_j), e_i\right>|\phi|^2\\
&=& {\rm Ric}^{V^\phi}_{ij},
\end{eqnarray*}
since $|\phi|=1$, where now $1\leq j\leq 2m$. 
On the other hand,
\begin{eqnarray*}
{\rm Re}\left<{\rm Ric}^{V^\phi}(e_j)\cdot\phi,e_i\cdot\phi\right> 
   &=&
{\rm Re}\left< \sum_{\alpha=1}^{2m}\sum_{1\leq k<l\leq r} \left<R^M(e_j ,e_\alpha)
e_{2m+k},e_{2m+l}\right>e_\alpha\cdot
\kappa_{r*}(f_kf_l)\cdot\phi,e_i\cdot\phi\right>\\
  &&
+ {\rm Re}\left<\sum_{\alpha=1}^{2m}\sum_{1\leq k<l\leq r} \Theta_{kl}(e_j ,e_\alpha)
e_\alpha\cdot\kappa_{r*}(f_kf_l)\cdot
\phi,e_i\cdot\phi\right>\\
   &&
  +{\rm Re}\left<{i}\sum_{\alpha=1}^{2m}\,dA(e_j ,e_\alpha)e_\alpha\cdot\phi,e_i\cdot\phi\right>\\
   &=&
\sum_{\alpha=1}^{2m}\,dA(e_j ,e_\alpha){\rm Re}\left<ie_\alpha\cdot\phi,e_i\cdot\phi\right>\\
   &=&
-\sum_{\alpha=1}^{2m}\,dA(e_j ,e_\alpha){\rm Re}\left<J^\phi(e_\alpha)\cdot\phi,e_i\cdot\phi\right>\\
   &=&
\sum_{\alpha=1}^{2m}\,dA(e_j ,e_\alpha)\left<J^\phi(e_\alpha),e_i\right>|\phi|^2\\
   &=&
\sum_{\alpha=1}^{2m}\,dA(e_j ,e_\alpha)J^\phi_{i\alpha}.
\end{eqnarray*}
Thus
\[{\rm Ric}^{V^\phi} =\left[{\rm proj}_{V^\phi}\circ\widehat{dA}|_{V^\phi}\right]\circ J^\phi,\]
where $\widehat{dA}$ denotes the skew-symmetric endomorphism  determined by $dA$ and metric
dualization.
\qd

{\bf Remark}. On each K\"ahler leaf, the spinor $\phi$ restricts to a parallel pure
Spin$^c$ spinor field.

\subsection{$D^{\perp}$-parallel partially pure spinor}\label{subsec: Dperp parallelness}

The following theorem is motivated by the condition
\[\nabla_u^S\phi =\lambda \, u \cdot \phi\]
for all $u\in\Gamma(D^\perp)$, $\lambda\in\mathbb{R}$, i.e. $\phi$ being a real $D^\perp$-Killing spinor.

\begin{theo}\label{theo: Dperp Killing}
Let $M$ be an oriented $n$-dimensional Riemannian manifold. The following are equivalent:
\begin{itemize}
\item[{\rm (i)}] $M$ admits a twisted Spin$^{c,r}(n)$ structure and a twisted spinor bundle $S$
% endowed with a connection $\nabla^S$ and
carrying a
partially pure spinor field $\phi\in\Gamma(S)$ satisfying
\[ (Y-iJ^\phi(Y))\cdot\nabla_{u}^S\phi=0\]
for every $Y \in \Gamma(V^\phi)$ and $u\in\Gamma((V^\phi)^\perp)$, where $\nabla^S$ is the covariant derivative described in subsection \ref{subsec: adapted connections}.
 \item[{\rm (ii)}] $M$ is endowed with an
almost CR hermitian structure of codimension $r$, where 
$D$ and $J$ are $D^\perp$-parallel. (In particular, the integral foliation of $D^\perp$ is totally
geodesic.) 
\end{itemize}
\end{theo}
{\em Proof.} 
Let us assume (i). For $X\in \Gamma(V^\phi)$,
\[X\cdot\phi=  i J^{\phi} X\cdot\phi\] 
Differentiate with respect to $u\in\Gamma((V^\phi)^\perp)$
\[\nabla_u X \cdot\phi + X \cdot \nabla_u^{S} \phi = i \nabla_u (J^\phi X)\cdot\phi 
+iJ^\phi X\cdot\nabla_u^{S}\phi,\]
so that
\[\nabla_u X \cdot\phi  = i \nabla_u (J^\phi X)\cdot\phi.\]
Since $\phi$ is a partially pure spinor 
\begin{eqnarray*}
 \nabla_u X&\in& V^\phi\\
\nabla_u (J^{\phi}X) &=& J^{\phi}(\nabla_u X),
\end{eqnarray*}
i.e. $D$ and $J$ are $D^\perp$ parallel, and so is $D^\perp$.

Conversely, if $M$ admits a CR hermitian structure, by Theorem \ref{theo: characterization almost CR}, 
$M$ admits a twisted Spin structure, a twisted spinor bundle $S$ endowed with a connection $\nabla^{S}$ as described in subsection \ref{subsec: adapted connections}, and a partially pure spinor
field $\phi\in\Gamma(S)$ such that $V^\phi=D$, $J^\phi=J$ 
and
\[(X-iJX)\cdot\phi =0\] 
for every $X\in\Gamma(D)$.
Let  $u\in\Gamma(D^\perp)$ and differentiate 
\[X\cdot\phi=  i JX\cdot\phi\] 
so that
\[\nabla_u X \cdot\phi + X \cdot \nabla_u^{S} \phi = i \nabla_u (JX)\cdot\phi 
+iJX\cdot\nabla_u^{S}\phi.\]
Since $J$ is $D^\perp$-parallel
\[\nabla_u (JX) = J(\nabla_uX),\]
and
\[ X \cdot \nabla_u^{S} \phi = iJX\cdot\nabla_u^{S}\phi.\]
i.e.
\[ (X- iJX)\cdot\nabla_u^{S}\phi=0.\]
\qd

{\bf Example}. The almost CR distribution $D_7$ on $\mathcal{G}_{m,1,r}$ gives the following example for Theorem \ref{theo: Dperp Killing}
\[
\begin{array}{ccc}
S^r & \hookrightarrow & {SO(2m+r+1)\over U(m)\times SO(r)}\\
 &  & \downarrow\\
 &  & {SO(2m+r+1)\over U(m)\times SO(r+1)},
\end{array}
\]

{\bf Remark}. 
 A generalized $D^\perp$-Killing partially pure spinor field $\phi$ is a spinor such that
\[\nabla_u^S\phi = E(u)\cdot \phi,\]
where $E$ is a symmetric endomorphism of $D^\perp$. Such a spinor also satisfies the hypotheses of Theorem \ref{theo: Dperp Killing}.

\vspace{.1in}

From Theorems \ref{theo: D-parallel} and \ref{theo: Dperp Killing} we obtain the following.

\begin{corol}
Let $M$ be an oriented $n$-dimensional Riemannian manifold. The following are equivalent:
\begin{itemize}
 \item[{\rm (i)}] $M$ is locally the Riemannian product of a K\"ahler manifold and a Riemannian manifold.
\item[{\rm (ii)}] $M$ admits a twisted Spin$^{c,r}(n)$ structure and a twisted spinor bundle
$S$ carrying a
partially pure spinor field $\phi\in\Gamma(S)$ satisfying
\[
(Y-iJ^\phi(Y))\cdot\nabla_{Z}^S\phi=0 
\]
for every $Y \in \Gamma(V^\phi)$ and $Z\in\Gamma(TM)$, where  $\nabla^S$ is the covariant derivative described in subsection \ref{subsec: adapted connections}.
\end{itemize}
\qd
\end{corol}

In the case of a real $D^\perp$-Killing partially pure spinor, we can say a little more about the foliation leaves' curvature.

\begin{theo}\label{Ricci-Killing-spinor}
Let $M$ be a Spin$^{c,r}$ $n$-dimensional Riemannian manifold such that its twisted spinor bundle $S$
admits a partially pure spinor field $\phi\in\Gamma(S)$ satisfying
\[ \nabla_{u}^S\phi=\mu\,u\cdot \phi\]
for every $u \in \Gamma((V^\phi)^\perp)$, where $\nabla^S$ is the connection described in subsection \ref{subsec: adapted connections} and
$\mu\in\mathbb{R}$, i.e. $\phi$ is real Killing in the 
directions of $(V^\phi)^\perp$.
Then 
\begin{itemize}
 \item The Ricci tensor decomposes as follows
\[
 {\rm Ric}^{(V^\phi)^\perp} = 4(r-1)\mu^2{\rm Id}_{(V^\phi)^\perp}+\sum_{1\leq k<l\leq r} 
\left[{\rm proj}_{(V^\phi)^\perp}\circ\hat{\Theta}_{kl}|_{(V\phi)^\perp}\right]\circ \hat\eta_{kl}^\phi,
\]
where $\Theta_{kl}$ denote the local curvature $2$-forms corresponding to the auxiliary connection on the $SO(r)$ principal bundle.
 \item The scalar curvature of each leaf tangent to $(V^\phi)^\perp$ is given by
\[{\rm R}^{(V^\phi)^\perp}
=4r(r-1)\mu^2+\sum_{1\leq k<l\leq r} 
\tr\left(\left[{\rm proj}_{(V^\phi)^\perp}\circ\hat{\Theta}_{kl}|_{(V\phi)^\perp}\right]\circ
\hat\eta_{kl}^\phi\right).
\]
 \item If
\[
  \sum_{1\leq k<l\leq r} 
\left[{\rm proj}_{(V^\phi)^\perp}\circ\hat{\Theta}_{kl}|_{(V\phi)^\perp}\right]\circ \hat\eta_{kl}^\phi
= \lambda \,{\rm Id}_{(V^\phi)^\perp}
\]
along a leaf of the foliation tangent to $(V^\phi)^\perp$ for some constant $\lambda\in\mathbb{R}$, then the
leaf is Einstein. 
\end{itemize}
\end{theo}
{\em Proof}.
Since $\phi$ is partially pure, $n=2m+r$ where ${\rm rank}(V^\phi)=2m$ and
${\rm rank}((V^\phi)^\perp)=r$.
Let $(e_1,\ldots, e_n)$ and $(f_1,\ldots f_r)$ be local orthonormal frames of $TM$ and $F$ respectively, such
that
\begin{eqnarray*}
V^\phi
   &=&
  \span(e_1,\ldots, e_{2m}),\\
e_{2j}
   &=& J^\phi (e_{2j-1}),\\
(V^\phi)^\perp
   &=&
  \span(e_{2m+1},\ldots, e_{2m+r}),\\ 
\eta_{kl}^\phi 
   &=&
  e_{2m+k}\wedge e_{2m+l},
\end{eqnarray*}
for $1\leq j\leq m$ and $1\leq k<l\leq r$.
First, if $u,v\in\Gamma((V^\phi)^\perp)$,
\begin{eqnarray*}
 R^{\theta,A}(u,v)\phi
   &=&
  \mu^2( v\cdot u
- u\cdot v)\cdot \phi
\end{eqnarray*}
Now, for $2m+1\leq i,j\leq 2m+r$,
\begin{eqnarray*}
 \sum_{i = 2m+1}^{2m+r} e_i\cdot R^{\theta,A}(e_j, e_i)(\phi) 
&=& 
 -2(r-1)\mu^2  e_j\cdot\phi.
\end{eqnarray*}
By taking the real part of the hermitian product with $e_t\cdot\phi$
we get
\begin{eqnarray*}
 {\rm Re}\left[-2(r-1)\mu^2 \left<e_j\cdot \phi,e_t\cdot \phi\right>\right] 
&=& -2(r-1)\mu^2 \delta_{jt}.
\end{eqnarray*}
If $u,v\in\Gamma((V^\phi)^\perp)$ then, by Theorem \ref{theo: Dperp Killing}, $[u,v]\in\Gamma((V^\phi)^\perp)$,
and
\begin{eqnarray*}
R^M(u,v)e_i 
   \in 
  \Gamma(V^\phi) \quad\quad&&\mbox{if $1\leq i \leq 2m$,} \\ 
R^M(u,v)e_i 
   \in 
  \Gamma((V^\phi)^\perp) \quad\quad&&\mbox{if $2m+1\leq i \leq 2m+r$.} 
\end{eqnarray*}
so that
\begin{eqnarray*}
\left<R^M(u,v)e_i, e_j\right>=0 
    \quad\quad&&\mbox{if $1\leq i \leq 2m$, $2m+1\leq j \leq 2m+r$,} \\ 
\left<R^M(u,v)e_i,e_j\right> =0
    \quad\quad&&\mbox{if $2m+1\leq i \leq 2m+r$, $1\leq j \leq 2m$.} 
\end{eqnarray*}
Now, if $1\leq\alpha\leq r$,
\begin{eqnarray*}
R^{\theta,A}(u,e_{2m+\alpha})\phi 
   &=& 
{1\over 2} \sum_{1\leq i<j\leq n} \left<R^M(u,e_{2m+\alpha})e_i,e_j\right> e_ie_j\cdot\phi\\
   &&
+{1\over 2} \sum_{1\leq k<l\leq r} \Theta_{kl}(u,e_{2m+\alpha}) \kappa_{r*}(f_kf_l)\cdot
\phi
  +{i\over 2}\,dA(u,e_{2m+\alpha})\phi\\
   &=& 
{1\over 2} \sum_{1\leq i<j\leq 2m} \left<R^M(u,e_{2m+\alpha})e_i,e_j\right> e_ie_j\cdot\phi\\
   &&
+{1\over 2} \sum_{1\leq k<l\leq r} \left<R^M(u,e_{2m+\alpha})e_{2m+k},e_{2m+l}\right>
e_{2m+k}e_{2m+l}\cdot\phi\\
  &&
+{1\over 2} \sum_{1\leq k<l\leq r} \Theta_{kl}(u,e_{2m+\alpha}) \kappa_{r*}(f_kf_l)\cdot
\phi
  +{i\over 2}\,dA(u,e_{2m+\alpha})\phi,
\end{eqnarray*}
where $dA$ denotes the curvature 2-form of the auxiliary connection on the $U(1)$-principal bundle.
Multiply by $e_{2m+\alpha}$ and sum over $\alpha$, $1\leq\alpha\leq r$,
\begin{eqnarray*}
\sum_{\alpha=1}^{r} e_{2m+\alpha}\cdot R^{\theta,A}(u,e_{2m+\alpha})\phi
   &=& 
\sum_{\alpha=1}^{r} \sum_{1\leq i<j\leq 2m} \left<R^M(u,e_{2m+\alpha})e_i,e_j\right>
e_{2m+\alpha}\cdot e_ie_j\cdot\phi\\
   &&
+\sum_{\alpha=1}^{r} \sum_{1\leq k<l\leq r} \left<R^M(u,e_{2m+\alpha})e_{2m+k},e_{2m+l}\right>
e_{2m+\alpha}\cdot e_{2m+k}e_{2m+l}\cdot\phi\\
  &&
+\sum_{\alpha=1}^{r} \sum_{1\leq k<l\leq r} \Theta_{kl}(u,e_{2m+\alpha})
e_{2m+\alpha}\cdot\kappa_{r*}(f_kf_l)\cdot
\phi\\
   &&
  +{i}\sum_{\alpha=1}^{r}dA(u,e_{2m+\alpha})e_{2m+\alpha}\cdot\phi.
\end{eqnarray*}
Furthermore,
\begin{eqnarray*}
&&{\rm Re}\left<\sum_{\alpha = 1}^r e_{2m+\alpha}\cdot R^{\theta,A}(e_{2m+\gamma},
e_{2m+\alpha})(\phi),e_{2m+\beta}\cdot\phi\right> \\
   &=& 
{\rm Re}\left<\sum_{\alpha=1}^{r} \sum_{1\leq i<j\leq 2m}
\left<R^M(e_{2m+\gamma},e_{2m+\alpha})e_i,e_j\right>
e_{2m+\alpha}\cdot e_ie_j\cdot\phi,e_{2m+\beta}\cdot\phi\right>\\
   &&
+{\rm Re}\left<\sum_{\alpha=1}^{r} \sum_{1\leq k<l\leq r}
\left<R^M(e_{2m+\gamma},e_{2m+\alpha})e_{2m+k},e_{2m+l}\right>
e_{2m+\alpha}\cdot e_{2m+k}e_{2m+l}\cdot\phi,e_{2m+\beta}\cdot\phi\right>\\
  &&
+{\rm Re}\left<\sum_{\alpha=1}^{r} \sum_{1\leq k<l\leq r} \Theta_{kl}(e_{2m+\gamma},e_{2m+\alpha})
e_{2m+\alpha}\cdot\kappa_{r*}(f_kf_l)\cdot
\phi,e_{2m+\beta}\cdot\phi\right>\\
   &&
  +{\rm
Re}\left<{i}\sum_{\alpha=1}^{r}dA(e_{2m+\gamma},e_{2m+\alpha})e_{2m+\alpha}\cdot\phi,e_{2m+\beta}
\cdot\phi\right>\\
   &=& 
-\left<{\rm Ric}^{(V^\phi)^\perp}(e_{2m+\alpha}), e_{2m+\beta}\right>|\phi|^2
\\
  &&
+\sum_{\alpha=1}^{r} \sum_{1\leq k<l\leq r} \Theta_{kl}(e_{2m+\gamma},e_{2m+\alpha})
{\rm Re}\left<e_{2m+\alpha}\cdot\kappa_{r*}(f_kf_l)\cdot
\phi,e_{2m+\beta}\cdot\phi\right>\\
   &=& 
-{\rm Ric}^{(V^\phi)^\perp}_{2m+\beta,2m+\gamma}
+\sum_{\alpha=1}^{r} \sum_{1\leq k<l\leq r} \Theta_{kl}(e_{2m+\gamma},e_{2m+\alpha})
\eta_{kl}^\phi(e_{2m+\alpha},e_{2m+\beta}),
\end{eqnarray*}
i.e.
\[
 {\rm Ric}^{(V^\phi)^\perp} = 4(r-1)\mu^2{\rm Id}_{(V^\phi)^\perp}+\sum_{1\leq k<l\leq r} 
\left[{\rm proj}_{(V^\phi)^\perp}\circ\hat{\Theta}_{kl}|_{(V\phi)^\perp}\right]\circ \hat\eta_{kl}^\phi.
\]
\qd

\subsection{CR foliation with Equidistant leaves}

\begin{theo}
Let $M$ be a Spin$^{c,r}$ Riemannian manifold such that its twisted spinor bundle $S$ admits a partially pure spinor field $\phi\in\Gamma(S)$.  
Let $V^\phi$ denote the almost-CR distribution and $(V^\phi)^{\perp}$ its orthogonal distribution.
If
\begin{eqnarray}
(X-iJ^{\phi}X)\cdot\nabla_u^S\psi &=& 0,\nonumber\\
(X-iJ^{\phi}X)\cdot \nabla_{X}^S\psi &=&0,\nonumber
\end{eqnarray}
for all $X\in \Gamma(V^\phi)$ and $u\in\Gamma((V^\phi)^{\perp})$,
where $\nabla^S$ the covariant derivative described in subsection \ref{subsec: adapted connections},
then
\begin{itemize}
\item $V^\phi$, $J^\phi$ and $(V^\phi)^\perp$ are $(V^\phi)^\perp$-parallel;
\item the totally geodesic foliation tangent to $(V^\phi)^\perp$ has equidistant leaves.
\end{itemize}
Furthermore, if the complex structure $J$ descends to the space of leaves $N$ at regular points, such a
complex structure is nearly-K\"ahler structure.
\end{theo}

{\em Proof}. The first statement follows from Theorem \ref{theo: Dperp Killing}.
Recall the condition for a foliation to have equidistant leaves 
\cite[Proposition 7]{Dearricott}
\begin{equation}
 \big< \nabla_X u , Y\big>+\big< X, \nabla_Y u \big>=0,\label{equidistant-leaves}
\end{equation}
for every $X,Y\in\Gamma(D)$ and $u\in \Gamma(D^{\perp})$.
Since $\big< u , Y\big>=\big< u , X\big>=0$,
\[\big< \nabla_X u , Y\big>+\big< u, \nabla_X Y \big>=0,\]
\[\big< \nabla_Y u , X\big>+\big< u, \nabla_Y X \big>=0,\]
so that \rf{equidistant-leaves} becomes
\[\big< \nabla_X Y + \nabla_Y X, u\big>=0.\]
We must prove $\nabla_X Y + \nabla_Y X\in \Gamma(D)$.
Taking covariant derivative with respect to $Y$ on
\begin{eqnarray}
 X\cdot \psi &=& i J(X)\cdot \psi,\nonumber
\end{eqnarray}
and with respect to $X$ on
\begin{eqnarray}
 Y\cdot \psi &=& i J(Y)\cdot \psi.\nonumber
\end{eqnarray}
we get
\begin{eqnarray}
 \nabla_Y X\cdot \psi+X\cdot \nabla_Y^S\psi &=& i \nabla_Y(J(X))\cdot \psi + i J(X)\cdot
\nabla_Y^S\psi,\nonumber\\
 \nabla_X Y\cdot \psi+Y\cdot \nabla_X^S\psi &=& i \nabla_X(J(Y))\cdot \psi + i J(Y)\cdot
\nabla_X^S\psi.\nonumber
\end{eqnarray}
Rearranging terms
\begin{eqnarray}
 \nabla_Y X\cdot \psi+(X - i J(X))\cdot
\nabla_Y^S\psi&=& i \nabla_Y(J(X))\cdot \psi ,\nonumber\\
 \nabla_X Y\cdot \psi+(Y- i J(Y))\cdot
\nabla_X^S\psi &=& i \nabla_X(J(Y))\cdot \psi .\nonumber
\end{eqnarray}
Adding up the last two equations and using 
\begin{eqnarray*}
 0
   &=& ((X+Y)-iJ(X+Y))\cdot \nabla_{X+Y}^S\psi\\ 
   &=& (X-iJ(X))\cdot \nabla_{X}\psi + (Y-iJ(Y))\cdot \nabla_{Y}^S\psi \\ 
   && + (X-iJ(X))\cdot \nabla_{Y}\psi + (Y-iJ(Y))\cdot \nabla_{X}^S\psi \\ 
   &=&  (X-iJ(X))\cdot \nabla_{Y}\psi + (Y-iJ(Y))\cdot \nabla_{X}^S\psi,
\end{eqnarray*}
we get
\[
 (\nabla_Y X+\nabla_X Y)\cdot \psi= i (\nabla_Y(J(X)) 
 +  \nabla_X(J(Y)))\cdot \psi .
\]
This means that $\nabla_Y X+\nabla_X Y\in\Gamma(D)$, i.e.
the foliation has equidistant leaves.
Furthermore, we have
\[
 \nabla_Y(J(X)) 
 +  \nabla_X(J(Y))
=  J(\nabla_Y X+\nabla_X Y).
\]
By setting $X=Y$, 
\[
 \nabla_X(J(X)) =  J(\nabla_X X) ,
\]
i.e.
\[(\nabla_XJ)(X)=0.\]

Since the leaves of the foliation are equidistant, by \cite[Lemma 19]{Dearricott}
the leaf space $N$ inherits a Riemannian metric at regular points and the quotient map $\pi$ is a
Riemannian submersion. 
The Levi-Civita connection at a regular point of $N$ is given by
\[\nabla^*_xy= \pi_*(\nabla_XY)\]
where $\pi_*(X)=x$ and $\pi_*(Y)=y$.

Finally, if we also have that $J$ descends to a complex structure $\mathcal{J}$ on $N$ 
in the form $\mathcal{J}(\pi_*(X)) = \pi_*(JX)$
and
\begin{eqnarray*}
 \nabla^*_x(\mathcal{J}(x)) -  \mathcal{J}(\nabla^*_x x)
   &=&  \pi_*\nabla_X({J}(X)) -  \mathcal{J}(\pi_*(\nabla_X X))\\ 
   &=&  \pi_*\nabla_X({J}(X)) -  \pi_*(J(\nabla_X X))\\ 
   &=&  \pi_*(\nabla_X({J}(X)) -  J(\nabla_X X))\\ 
   &=&  \pi_*(0)\\ 
   &=&  0. 
\end{eqnarray*}
i.e. $\mathcal{J}$ is a nearly-K\"ahler structure at regular points of $N$.
\qd

\section{Codimension 1 almost CR structures}\label{sec: codimension one CR structures}

Now we will focus our attention on oriented Riemannian manifolds of dimension $n = 2m+1$, i.e. $r=1$.
Recall that
\[Spin^{c,1}(n)=Spin^c(n),\]
so we are back to Spin$^c$ geometry.

Let $M$ be a $2m+1$-dimensional Spin$^c$ manifold carrying a partially pure spinor field $\psi$.
We have 
\[TM = V^\psi\oplus (V^\psi)^\perp,\] 
where $(V^\psi)^\perp$ is a trivial real line bundle over $M$. 
Consider the vector field $\xi^\psi$ defined by 
\[\left<X, \xi^\psi\right> = i \left<X\cdot\psi, \psi\right>,\]
for all $X \in \Gamma(TM)$. 
If $X \in \Gamma(V^\psi)$,  
\[X\cdot\psi=  iJ^\psi(X)\cdot\psi,\] 
so that taking the scalar product with $\psi$ gives
\[\left<X\cdot\psi, \psi\right> = i \left<J^\psi(X)\cdot\psi, \psi\right> = 0,\]
i.e. 
\begin{equation}
\left<X, \xi^\psi\right> = 0.\label{eq: Vpsi}
\end{equation} 
Let $\theta^\psi=(\xi^\psi)^*$ the metric dual $1$-form, and $G_\psi$ the symmetric $2$-form defined by 
\[G_\psi(X, Y) = d\theta^\psi (X, J^\psi(Y)),\]
for any $X, Y \in \Gamma(V^\psi)$. 

\subsection{Pseudoconvex CR manifolds}

\begin{theo} Let $M$ be an oriented $2m+1$-dimensional smooth manifold. Then, 
 $M$ is a pseudoconvex CR manifold if and only if it has a Spin$^c$ structure carrying 
 an integrable partially pure spinor such that  $G_\psi$  is
positive definite.
\label{pseudoCR}
\end{theo}
{\em Proof.} Assume that $M$ is a pseudoconvex CR manifold. Then $D = \ker\ \theta$ for some
hermitian structure $\theta$. We consider the Tanaka-Webster (Riemannian) metric $g_{\theta}$. 
It is known that  $(M, g_{\theta})$ has a canonical Spin$^c$ structure
carrying an integrable partially pure spinor field $\psi$ of unit length for which $D=V^\psi$. 
Moreover, we recall that there exists a unique vector field $T$ such that $\theta (T)
=1$ and $T\lrcorner d\theta = 0$. 
We claim that $T = \xi^\psi$ so that $\theta = \theta^\psi$ and $G_\psi$ is positive definite. 

Let us note that $\xi^\psi\not =0$, since
\[T \cdot\psi = -i\psi\] 
implies
\begin{eqnarray}
g_\theta(T,\xi^\psi)
   &=&
i\left<T\cdot \psi,\psi\right>\nonumber\\
   &=& 
  1.\label{eq: <T,xi>}
\end{eqnarray}
Since $g_\theta(X, \xi^\psi) = 0$ for all $X \in\Gamma(D)$,
$\xi^\psi $ is a multiple of $T$
and by \rf{eq: <T,xi>}
\[T = \xi^\psi.\] 

Conversely, assume that a $2m+1$ dimensional Riemannian manifold
$M$ carries an integrable partially pure spinor field and
such that $G_\psi$ is positive definite. 
It remains to prove that it is a pseudoconvex structure. For
this, it is sufficient to prove that $\ker \theta^\psi = V^\psi$
which is already the case by \rf{eq: Vpsi}.
\qd

{\bf Remark}.
From the proof of Theorem \ref{pseudoCR}, 
we see that the
Tanaka Webster metric   is given by  
\[g_{\theta^\psi} (X, Y) = G_\psi (X, Y)\]
for $X,Y\in V^\psi$.

\subsection{Metric contact and Sasakian manifolds}

Let $M^{2m+1}$ be an oriented smooth manifold and
$(\mathfrak{X}, \xi, \eta)$
a synthetic object consisting of a $(1,1)$-tensor field $\mathfrak{X}: TM \longrightarrow TM$, a tangent
vector field $\xi$, and a differential  $1$-form $\eta$ on $M$. $(\mathfrak{X}, \xi, \eta)$ is an {\em almost
contact struture} if
$$\mathfrak{X}^2= -\mathrm{Id}+\eta\otimes\xi,\ \ \ \mathfrak{X}\xi=0, \ \ \eta(\xi)=1,\ \ \eta\circ\mathfrak{X}=0.$$ 

An almost contact structure is said to be {\em normal} if, for all $X, Y\in \Gamma(TM)$, $N^\mathfrak{X}  =  d\eta\otimes \xi,$
where $N^\mathfrak{X}$ is the Nijenhuis contact tensor of $\mathfrak{X}$ defined by
$$N^\mathfrak{X} (X, Y) = -\mathfrak{X}^2[X, Y] + \mathfrak{X}[\mathfrak{X} X, Y] + \mathfrak{X}[X, \mathfrak{X} Y] - [\mathfrak{X} X, \mathfrak{X} Y],$$ 
for all $X, Y\in \Gamma(TM)$. 

A Riemannian metric $g$ is said to be compatible with the almost
contact structure if
$$g(\mathfrak{X} X,\mathfrak{X} Y) = g(X, Y) -\eta(X)\eta(Y),$$
for all $X, Y\in \Gamma(TM)$. An almost contact structure $(\mathfrak{X}, \xi, \eta)$ together with a
compatible
Riemannian metric $g$ is called an {\em almost contact metric structure}. 

Given an an almost contact
metric structure $(\mathfrak{X}, \xi, \eta, g)$ one defines a $2$-form $\alpha$ by $\alpha(X, Y) = g(X, \mathfrak{X}
Y)$ for all $X, Y\in \Gamma(TM)$. Now, $(\mathfrak{X}, \xi, \eta, g)$ is said to satisfy the {\em contact
condition} if $\alpha=d\eta$ and if it is the case, $(\mathfrak{X}, \xi, \eta, g)$ is called a {\em contact metric
structure} on $M$. 

A contact metric structure $(\mathfrak{X}, \xi, \eta, g)$ that is also normal is called a
{\em Sasakian structure} (and $M$ is {\em Sasakian}). 

Every strictly pseudoconvex CR manifold is also a contact metric manifold \cite{dragomir, blair}.
Moreover, this contact metric  structure is a Sasakian structure if and only if  the Tanaka torsion
vanishes, i.e., $\tau =0$ \cite{dragomir, blair}. Conversely, a metric contact manifold has a
natural almost CR structure \cite{dragomir, blair}. This almost CR structure is a CR structure (and
then automatically strictly pseudoconvex) if 
and only if  $\mathfrak{X} \circ N^\mathfrak{X} (X, Y)= 0$ for all $X, Y \in \Gamma(\ker\eta)$ \cite{nic}.

\begin{corol}
Every oriented contact Riemannian manifold has a Spin$^c$ structure carrying a partially pure spinor.
Moreover, 
this spinor is integrable if  the contact structure is normal, i.e., if $M$ is Sasakian. 
\end{corol}
{\em Proof.} We know that any contact Riemannian manifold has a Spin$^c$ structure. Moreover, for
this Spin$^c$ structure carries a partially pure
spinor. We have
\begin{eqnarray*}
(N^\mathfrak{X} (X, Y) +\eta([X, Y]) \xi)\cdot \psi = (X-i\mathfrak{X} X)\nabla_{Y-i\mathfrak{X} Y} \psi - (Y-i\mathfrak{X}
Y)\nabla_{X-i\mathfrak{X}
X} \psi,
\end{eqnarray*}
for all $X, Y\in \Gamma(V^\psi)$. But $d \eta(X, Y) = -\eta([X, Y])$ for all $X, Y \in
\Gamma(\mathrm{ker}\eta)$. Hence, if the contact metric is normal, $\psi$ is integrable.\qd

\begin{corol}
Let $M$ be an oriented contact Riemannian manifold.
Then, it is a Spin$^c$ manifold carrying an integrable 
partially pure spinor $\psi$ with positive definite $G_\psi$ if and
only if $\mathfrak{X}\circ N^{\mathfrak{X}}=0$.
\qd
\end{corol}

\begin{corol}
 Let $M$ be a Sasakian manifold satisfying  $N^{\mathfrak{X}}\circ \mathfrak{X}=0$. 
 Then, it is a Spin$^c$ manifold
carrying a integrable partially pure spinor field with positive definite $G_\psi$. 
Conversely, if $M$ is a
Riemannian Spin$^c$ manifold carrying a integrable 
strictly partially pure spinor field with positive definite $G_\psi$ and  
such that  $\tau = 0$, then $M$ is a Sasakian manifold.
\qd
\end{corol}

\subsection{Isometric immersions via partially pure Spin$^c$ spinors}\label{6}

Let $N^{2m-1}$ be an oriented real hypersurface of a K\"{a}hler manifold $(M^{2m}, \overline g,
J)$ endowed with the metric $g$ induced by $\overline g$.
We denote by $\nu$ the unit normal inner vector globally defined on $M$ and by ${\rm II}$ the second
fundamental form of the 
immersion. 

The complex structure $J$ induces an almost contact metric structure 
$(\mathfrak{X},\xi,\eta, g)$ on $M$, 
where $\mathfrak{X}$ is the $(1,1)$-tensor defined by 
\[g(\mathfrak{X} X,Y)=\overline{g}(JX,Y)\] 
for all
$X,Y\in\Gamma(TN)$, $\xi=-J\nu$ is a tangent vector field and $\eta$ is
the $1$-form associated with $\xi$, that is  
\[\eta(X)=g(\xi,X).\] 

For every $X\in \Gamma(TN)$
\begin{eqnarray*}
\mathfrak{X}^2X=-X+\eta(X)\xi,\quad g(\xi,\xi)=1,\quad\text{and}\quad \mathfrak{X}\xi=0.
\end{eqnarray*}
Moreover, from the relation between the Riemannian connections $\overline\nabla$ of $M$ and
$\nabla$ of $N$, 
\[\overline{\nabla}_XY=\nabla_XY+g({\rm II}(X),Y)\nu,\] we deduce the following two
identities:
\begin{eqnarray*}
(\nabla_X\mathfrak{X})Y&=&\eta(Y){\rm II}(X)-g({\rm II}(X),Y)\xi \\
\nabla_X\xi&=&\mathfrak{X}( {\rm II}(X)), 
\end{eqnarray*}
for every $X, Y \in \Gamma(TN)$. 
We can choose an orthonormal frame 
\[\mathcal{B}=\{e_1, e_2 = \mathfrak{X}
e_1, \dots, e_{2m-3}, e_{2m-2} = \mathfrak{X} e_{2m-3}, \xi\}\] 
of $N$ such that  
\[\mathcal{B}\cup \{\nu = J\xi\}\]
is an orthonormal frame of $M$.

\begin{theo}
Let $(M^{2m}, g, J)$ be a Hermitian manifold. Then, any real oriented hypersurface $N\subset M$ has a
Spin$^c$ structure carrying an integrable 
partially pure spinor. 
\label{kaa}
\end{theo}
{\em Proof}. Since $M$ is a Hermitian manifold, it has a canonical Spin$^c$ structure carrying a 
integrable pure spinor field $\psi$. The restriction of this Spin$^c$ structure to
an oriented real hypersurface $N^{2m-1}$ gives a Spin$^c$ structure carrying a spinor field
$\phi=\psi\vert_N$ satisfying
\[\nabla^N_X \phi = \nabla^M_X \psi \vert_M -\frac 12 {\rm II}(X)\bullet\phi\]
for all $X \in \Gamma(TN)$.
We claim that the spinor field $\phi$ is an integrable partially pure spinor field.  
For $j=1, \cdots, 2m-2$, we have 
\begin{eqnarray*}
(e_j - i \mathfrak{X} e_j)\bullet\phi 
  &=& 
   (e_j - i \mathfrak{X} e_j)\cdot \nu\cdot\psi\vert_M \\
  &=& 
   - \nu \cdot (e_j - i \mathfrak{X} e_j)\cdot\psi\vert_M \\
  &=& 
   0, 
\end{eqnarray*}
since 
\begin{eqnarray*}
(e_j - i \mathfrak{X} e_j) \cdot\psi 
   &=& 
    (e_j - i J e_j)\cdot\psi \\
   &=&
    0.
\end{eqnarray*}
The distribution 
\[V^\phi =\{ X \in \Gamma(TN), \ \ X\bullet\phi=  - i \mathfrak{X} X \bullet\phi\},\]
is of constant rank $2(m-1)$
and  is an almost CR structure 
By the Spin$^c$ Gauss formula, the spinor $\phi$ is
integrable if and only if $\psi$ is integrable because 
\[(X-iJX)\bullet {\rm II} (Y- iJY)\bullet \phi = (Y-iJY)\bullet {\rm II} (X- iJX)\bullet\phi,\]
for all $X, Y\in \Gamma(V^\phi)$. Since $\psi$ is integrable, $\phi$ is also
integrable.\qd

{\bf Remark}.
From Proposition \ref{kaa}, any real oriented hypersurface $N$ of a  K\"{a}hler manifold $M$ has a
Spin$^c$ structure carrying an integrable partially pure spinor $\phi$
satisfying
$$\nabla^N_X \phi = -\frac 12 {\rm II}(X)\bullet\phi,$$
for all $X\in \Gamma(TN)$, i.e. a generalized Killing spinor.

\begin{theo}\label{immersion}
Let $(N^{2m+1}, g)$ be an oriented almost contact metric  Spin$^c$ manifold  carrying a
parallel partially pure
 spinor $\psi$, and $I=[0,1]$. Then the product $\mathcal{Z} :=M\times I$
endowed with the metric $dt^2+g$ and the Spin$^c$ structure arising from the given one on
$M$ is a K\"{a}hler manifold having a parallel spinor $\psi$ whose restriction to $M$ is just
$\phi$.
\end{theo}
{\em  Proof}.  
First, the pull back of the Spin$^c$ structure on $M$ defines a Spin$^c$ structure on
$M\times I$. Moreover, from the spinor field $\phi$, we can construct on $M\times I$ a parallel
spinor $\psi$. It remains to show that $M \times I$ is K\"ahler. We define the endomorphism
$\overline J$ by 
\begin{eqnarray*}
\overline J(X) &=& J(X)\quad\text{for any}\quad X\in\Gamma(D), \\  
\overline J(T)&=&\nu, \\  
\overline J(\nu) &=& -T. 
\end{eqnarray*}
It is easy to prove that $(M\times I, \overline J, g+dt^2)$ is an almost Hermitian manifold.
Moreover, since $T\bullet\phi = - i\phi$, then 
$$\nabla_T T = \nabla_X T = \nabla_{JX}T = 0.$$
Since the immersion is totally geodesic, we get
\begin{eqnarray}
\label{lc1}
\overline\nabla_T T = \overline\nabla_X T = \overline\nabla_{JX}T = 0.
\end{eqnarray}
Now, since $\phi$ is parallel on $M$, we get for any $X\in \Gamma(D)$,
\begin{eqnarray}\label{lc2}
 \nabla_X J &=& 0,  \nonumber\\   
 \nabla_T X &\in& \Gamma(D), \label{lc2} \\
 \nabla_T JX &=& J(\nabla_TX).\nonumber
\end{eqnarray}
Finally, using (\ref{lc1}) and (\ref{lc2}), we conclude that $\overline{\nabla}\, \overline{J} = 0$ on
$M\times I$.\qd

{\small
\renewcommand{\baselinestretch}{0.5}
\newcommand{\bi}{\vspace{-.05in}\bibitem} }


\begin{thebibliography}{10}
\newcommand{\BHNY}{Berlin\kern3pt Heidelberg\kern3pt New York}

\bibitem{blair}  Blair, D. E.:{\it Contact manifolds in Riemannian geometry}, vol. 509, Springer-Verlag, Berlin-Heidelberg-New York, 1976.
\bibitem{borisov} Borisov, L.; Salamon, S.;  Viaclovsky, J.: {\it Twistor geometry and warped product
orthogonal complex structures},  Duke Math. J. 156 (2011), No. 1, 125-166.
\bibitem{cartan} Cartan, E.: {\it The theory of spinors}, Hermann, Paris, 1966.
\bibitem{cartan1} Cartan, E.: {\it Lecons sur la th\'eorie des spineurs}, Paris, Hermann (1937).
\bibitem{cheva}  Chevalley, C.: {\it The Algebraic Theory of Spinors}, Colombia University Press, New York, 1954.
\bibitem {dadok} Dadok, J.;  Harvey, R.:{\it Calibrations on $\mathbb{R}^6$}, Duke Math. J. 4 (1983),
1231-1243.
\bibitem{Dearricott} 
Dearricott, O.: {\em Lectures on n-Sasakian manifolds}. Geometry of manifolds with non-negative sectional curvature, 57-109, Lecture Notes in Math., {\bf 2110}, Springer, Cham, 2014
\bibitem{dragomir} Dragomir, S.; Tomassini, G.: {\it Differential geometry and analysis on CR manifolds}, Progress in mathematics, volume 246.

\bibitem{Espinosa-Herrera}  Espinosa, M.; Herrera, R.: {\it Spinorially twisted Spin structures, I: curvature identities
and eigenvalue estimates}. Differential Geom. Appl. {\bf 46} (2016), 79-107. 

\bibitem{Friedrich} Friedrich, T.:{\it Dirac operator's in Riemannian geometry},
Graduate studies in mathematics, Volume 25, American Mathematical Society.
\bibitem{har}  Harvey, R.; Lawson, H. B.: {\it Calibrated geometries}, Acta Math. 148 (1982),
47-157.

\bibitem{Herrera-Tellez} Herrera, R.; Tellez, I.: {\em Twisted partially pure spinors}. 
J. Geom. Phys. {\bf 106} (2016), 6-25
\bibitem{hitchin} Hitchin, N.: {\it Harmonic spinors}, Adv. Math. 14 (1974) 1-55.
\bibitem{morel} Morel, B.: {\it Tenseur d'impulsion-\'energie et g\'eom\'etrie spinorielle extrins\`eque},
Ph.D. thesis, Institut Elie Cartan, 2002.
\bibitem{morgan} Morgan, J. W.: {\em The Seiberg–Witten equations and applications to the topology of
smooth four-manifolds}. Mathematical Notes, vol. 44,
Princeton University Press, Princeton, NJ, ISBN: 0-691-02597-5, 1996, viii+128 pp.
\bibitem{r2} Nakad, R.: {\it The Energy-Momentum tensor on $\Spinc$ manifolds},  IJGMMP  Vol. 8, No. 2, 2011.
\bibitem{nic} Nicolaescu, L.: {\it Geometric connections and geometric Dirac operators on contact
manifolds.} Differential Geometry and its Applications, 22, 355-378, 2005.
\bibitem{pen1} Penrose, R.: {\it Twistor theory, its aims and achievements}, Quantum Gravity: an
Oxford Symposium (eds. C. J. Isham, R. Penrose, and D. Sciama), pages 268-407. Oxford University
Press, Oxford, 1975.
\bibitem{pen2} Penrose, P.; Rindler, W.: {\it Spinors and space-time}, volume 1. Cambridge University Press, Cambridge, 1986.
\bibitem{pen3} Penrose, R.; Rindler, W.: {\it Spinors and space-time}, volume 2. Cambridge University Press, Cambridge, 1986.
\bibitem{petit} Petit, R.: {\it Spin$^c$ structures and Dirac operators on contact manifolds}, Differential Geometry and its Applications 22 (2005), 229-252.
\bibitem{witten} Witten, E.: {\it Monopoles and four-manifolds}, Math. Res. Lett. 1 (6) (1994)
769–796.


\end{thebibliography}
\end{document}